\documentclass[11pt]{article}
\usepackage[utf8]{inputenc}
\usepackage{graphicx}
\usepackage{amsfonts,amssymb,amsmath,amsthm}
\usepackage{latexsym}
\usepackage{epstopdf}
\usepackage{mathrsfs}
\usepackage[all]{xy}
\usepackage{float} 
\usepackage[pagebackref]{hyperref}

\usepackage{color}
\definecolor{lav}{rgb}{0.80,0.0,0.95} 
\definecolor{cobre}{rgb}{.90,0.38,0.09} 
\definecolor{gren}{rgb}{0.28,.90,0.48} 
\def\C{{\mathbb C}}
\def\R{{\mathbb R}}
\def\s{\mathbb{S}} 
\def\cS{\mathcal{S}} 
\def\B{{\mathbb B}}
\def\D{{\mathbb D}}
\def\e{\varepsilon}
\def\a{\alpha}      
\def\0{{\underline 0}}

\def\codim{\operatorname{codim}}
\def\rank{\operatorname{rank}}
\def\Sing{\operatorname{Sing}}

\theoremstyle{plain}
\newtheorem{theorem}{Theorem}[section]
\newtheorem{corollary}[theorem]{Corollary}
\newtheorem{proposition}[theorem]{Proposition}
\newtheorem{lemma}[theorem]{Lemma}
\newtheorem{notation}[theorem]{Notation}

\theoremstyle{definition}
\newtheorem{definition}[theorem]{Definition}
\newtheorem{remark}[theorem]{Remark}
\newtheorem{example}[theorem]{Example}
\newtheorem{convention}[theorem]{Convention}

\def\rond{\mathaccent"7017}  

\numberwithin{equation}{section}

\title{The topology of complex map-germs and \\ a general L\^e-Greuel formula}
\author{L\^e D\~ung Tr\'ang\footnote{Professor L\^e D\~ung Tr\'ang passed away in November 2025. The present article had been essentially completed only two weeks before his passing.}, Juan J.   Nu\~no-Ballesteros, Jos\'e Seade}
\date{}

\begin{document}

\maketitle

\begin{abstract} 
Let $(X,\0)\subset(\mathbb C^N,\0)$ be a singular complex analytic germ
and let $f:(X,\0)\to(\mathbb C,0)$ be holomorphic. Given a holomorphic
germ $g:(X,\0)\to(\mathbb C,0)$ having an isolated critical point
relative to $f$, we use stratified Morse theory and relative polar
geometry to describe the topology of the Milnor fiber $F_f$ in terms
of the lower-dimensional fiber $F_{g,f}$ and the Morse data of a
morsification of $g|_{F_f}$. More precisely, we show that $F_f$ is
obtained from $F_{g,f}\times D^2$ by attaching stratified Morse data
whose multiplicities are determined by the relative polar geometry.

Taking Euler characteristics yields a stratified L\^e--Greuel formula
in which the contribution of each stratum combines a relative polar
multiplicity with the topology of its complex link. When $X$ has an
isolated singularity, this formula admits an algebraic interpretation
in terms of a finite-dimensional algebra defined by Jacobian ideals.
The resulting framework applies beyond complete intersections and,
notably, does not require $f$ to have an isolated critical point.
It also recovers several known L\^e--Greuel type formulas as special
cases.
\end{abstract}



\section*{Introduction} 

Understanding the topology of the fibers of a complex analytic function near a singular point is a fundamental problem in complex geometry and topology. One of the main approaches to this problem is to study the Milnor fiber by means of lower-dimensional sections and to determine how its topology can be reconstructed from theirs. This viewpoint has its roots in the ideas of S. Lefschetz and R. Thom, and was substantially developed by Teissier, L\^e, Hamm, and others in the 1970s. Polar geometry and vanishing cycles provide a natural framework for this comparison, and, in the classical complete intersection setting, the L\^e--Greuel formula gives a remarkable numerical relation between the corresponding topological invariants.

The aim of this paper is to develop this picture for Milnor fibers of
holomorphic functions on singular analytic spaces, without assuming
that the ambient space is a complete intersection or that the function
has isolated critical locus. Our results first establish the geometric
and stratified structure needed for this comparison and then give a
topological description of the Milnor fiber in terms of a
lower-dimensional slice and stratified Morse data. The resulting decomposition leads naturally to a general L\^e--Greuel type formula.

Let $(X,\0)\subset(\mathbb C^N,\0)$ be a reduced complex analytic germ,
endowed with a complex analytic Whitney stratification \cite{W}, and let
$f:(X,\0)\longrightarrow(\mathbb C,0)$
be holomorphic. We allow both $X$ and the critical locus of $f$ to be
singular and non-isolated. Let
$g:(X,\0)\longrightarrow(\mathbb C,0)$
be another holomorphic germ having an isolated critical point relative
to $f$. We denote by $F_f$ the Milnor fiber of $f$ and by $F_{g,f}$
the Milnor fiber of the map
\[
\Phi=(g,f):(X,\0)\longrightarrow(\mathbb C^2,\0).
\]

A first step is to understand the geometry of $\Phi$ under the sole
assumption that $g$ has an isolated critical point relative to $f$.
Although this condition is formulated in terms of critical points on
the fibers of $f$, it has strong consequences for the geometry of
$\Phi$. We prove that the relative critical locus away from
$f^{-1}(0)$ is either empty or a curve, that its image by $\Phi$ is an
analytic curve, and, more importantly, that $\Phi$ admits a suitable
stratification satisfying the Thom condition. This provides the
geometric framework needed to compare the Milnor fiber $F_f$ with the
lower-dimensional fiber $F_{g,f}$.

After replacing $g$ by a suitable morsification $g_t$, the restriction of $|g_t|$ to the fiber $F_f$
is a non-depraved stratified Morse function whose critical points are
controlled by the relative polar geometry of the pair $(g,f)$. For
each stratum $S_\alpha$ of the induced stratification of $F_f$, let
$\delta_\alpha$ denote the number of Morse critical points of $g_t$
on $S_\alpha\cap F_f$.

Our main topological theorem states that $F_f$ is obtained from
$F_{g,f}\times D^2$
by attaching, for each stratum $S_\alpha$, $\delta_\alpha$ copies of
the corresponding stratified Morse data. In particular, up to
homotopy the normal part of each attachment is determined by the
complex link of the stratum. Thus the topology of $F_f$ is governed
by two complementary ingredients: the multiplicities arising from
relative polar geometry and the topology of the complex links of the
strata.

This approach goes back to the seminal work of L\^e D\~ung Tr\'ang \cite{L1,L3} in the smooth ambient setting.
In a different direction, Massey's theory of L\^e cycles and L\^e numbers and his sheaf-theoretic study of Milnor fibers \cite{Massey0, Massey} provide powerful tools for describing the topology and cohomology of Milnor fibers, including in the presence of non-isolated singularities. For functions with an isolated stratified singularity on a singular ambient space, Tib\u{a}r's Bouquet Theorem \cite{Tibar} gives a description of the homotopy type of the Milnor fiber based on an underlying handlebody construction. 

Our result is closely related to this circle of ideas, but differs both in its hypotheses and in the topological information it retains. We allow the critical locus of $f$ to be non-isolated and take $g$ to be an arbitrary holomorphic germ having an isolated critical point relative to $f$; moreover, before passing to homotopy, we reconstruct the homeomorphism type of $F_f$ from $F_{g,f}\times D^2$ by attaching the full stratified Morse data, with multiplicities determined by the relative polar geometry.

The topological decomposition above has an immediate numerical consequence. Namely,
\[
\chi(F_f)-\chi(F_{g,f})
   =\sum_\alpha \delta_\alpha m_\alpha,
\]
where
\[
m_\alpha
   =(-1)^{n_\alpha}
     \chi\bigl(\operatorname{Cone}(\mathcal L_\alpha),\mathcal L_\alpha\bigr),
\]
and $\mathcal L_\alpha$ denotes the complex link of the stratum
$S_\alpha\cap F_f$. This gives a stratified L\^e--Greuel formula in
which the contribution of each stratum separates naturally into a
polar multiplicity $\delta_\alpha$ and a topological factor determined
by its complex link.
 In this way, the topological decomposition provides the underlying geometric picture, while the L\^e--Greuel formula captures its Euler-characteristic content.

When $\Sing(X)\subseteq f^{-1}(0)$ (and hence, when $(X,\0)$ has an isolated singularity), this formula admits a further
algebraic interpretation. If $\dim X=n+1$, we prove that
\[
\chi(F_f)-\chi(F_{g,f})
 =
(-1)^n
\dim_{\mathbb C}
\frac{\mathscr O_{X,\0}}
     {(f)+J_X(f,g):f^\infty}.
\]
Thus the topological difference between the Milnor fiber and its slice
is measured by the dimension of a naturally defined finite-dimensional
algebra. This provides a direct bridge between the stratified
Morse-theoretic description of the Milnor fiber and local analytic
invariants.

The scope of the results in this article goes beyond the classical setting in two directions: the ambient space need not be a complete intersection, and the function $f$ need not have an isolated critical point. 

Several previously known L\^e--Greuel type formulas arise naturally as special cases. In the setting of functions on singular spaces, related formulas were obtained by Dutertre and Grulha
\cite{DN}. Our results extend this circle of ideas in a different direction: the stratified Morse-theoretic reconstruction of the Milnor fiber leads, when $X$ has an isolated singularity, to an algebraic L\^e--Greuel formula that remains valid when $f$ has non-isolated critical locus, thereby extending the classical formula beyond the isolated-critical-point setting. This algebraic formulation also provides a unified framework for several applications, including isolated determinantal singularities, smoothable curve and surface singularities, and formulas related to image Milnor numbers. In these settings the topology of the relevant smoothing or Milnor fiber need not be concentrated in a single dimension, so that the usual description in terms of a number of vanishing spheres is no longer available.

From this point of view, the results of the paper reveal a common
geometric mechanism:
\[
\text{Milnor fiber}
\quad = \quad
\text{slice}
\quad + \quad
\text{stratified Morse data}.
\]
The resulting topological decomposition and the L\^e--Greuel formula
are two complementary manifestations of this mechanism. The former
describes how the Milnor fiber is assembled from a lower-dimensional
slice, while the latter measures this passage through the contributions
of the different strata, combining relative polar multiplicities with
the topology of their complex links. 

This framework also places several existing L\^e--Greuel type formulas in a common geometric setting; see, for instance,
 \cite{CMSS, DN, Massey, NOT, NP, Zach}.

The paper is organized as follows. Section~1 contains background
material on complex analytic functions and stratified Morse theory.
Section~2 discusses relative polar sets. Sections~3 and~4 develop the
geometric and topological ingredients needed for the main results.
In Section~5 we establish the structure of the relative critical locus
and discriminant, the required Thom stratification, and the topological
description of the Milnor fiber in terms of the slice and the
corresponding stratified Morse data. Section~6 gives the algebraic
interpretation of the resulting L\^e--Greuel formula in the case of
isolated singularities. Finally, Section~7 discusses applications and
recovers several known formulas as special cases.

 \subsection{Acknowledgements}
  
  Part of this work was completed while the first and third authors were visiting Universidad de Valencia, whose support and hospitality they gratefully acknowledge. The authors also thank Bernard Teissier,  Helmut Hamm and Roberto Gim\'enez Conejero for helpful conversations.

 J.J. Nu\~no-Ballesteros has been partially supported by Grant PID2025-172404NB-I00 funded by MCIN/AEI/10.13039/501100011033 and by ``ERDF A way of making Europe" and by Grant UV-INV\_AE\_4225594.
  Jos\'e Seade received support from UNAM-PAPIIT project IN101424.

\section{Preliminaries on singular points of complex analytic functions}\label{S: singularities}

Throughout this work, we denote by $\0$ the origin in the source
$X \subset \C^N$ and by $0$ the origin in the target. Similarly, we use $\B$ for balls in the source and $\D$ for balls in the target.

\subsection{Complex analytic functions on a singular complex analytic set}

Milnor's foundational work \cite{M} studies the local topology of
holomorphic functions with isolated critical points. Let $U$ be a
neighborhood of $\0\in\C^{n+1}$ and let
$
f:(U,\0)\to(\C,0)
$
be a holomorphic germ with an isolated critical point at $\0$. Then,
for sufficiently small $\varepsilon>0$, the sphere
$\mathbb S_\varepsilon$ meets $f^{-1}(0)$ transversely, and, for
$\eta>0$ sufficiently small, the restriction
$
f:\B_\varepsilon\cap f^{-1}(\D_\eta\setminus\{0\})
\to \D_\eta\setminus\{0\}
$
is a locally trivial fibration. Its fiber $F_f$, called the
\emph{Milnor fiber}, has the homotopy type of a bouquet of $\mu(f)$
spheres of real dimension $n$.

We now consider holomorphic functions on singular analytic spaces,
following \cite{L,L1,GM}. Let $(X,\0)$ be the germ of a closed complex
analytic subset of an open set $U\subset\C^N$, and let
$
f:(X,\0)\to(\C,0)
$
be a holomorphic germ. Fix a Whitney stratification
$
\mathcal S=\{S_\alpha\}_{\alpha\in A}
$
of $X$.

\begin{definition}\label{def: singular point}
A point $x\in X$ is \emph{regular} for $f$ if the restriction of $f$
to the stratum containing $x$ is a submersion at $x$. Otherwise, $x$
is a \emph{critical point}. A critical point $x$ is \emph{isolated}
if it has a neighborhood containing no other critical points.
\end{definition}

To study such points, we recall the conormal space of $X$:
\begin{equation}\label{d. conormal}
T_X^*U:=\bigcup_{\alpha\in A}T^*_{S_\alpha}U,
\end{equation}
where
$
T^*_{S_\alpha}U
=
\{(x,\xi)\in T^*U\mid x\in S_\alpha,\ \xi(T_xS_\alpha)=0\}.
$
By Whitney condition (a), $T_X^*U$ is a closed analytic subset of
$T^*U$ of dimension $N$, naturally stratified by the pieces
$T^*_{S_\alpha}U$.

We regard $d\widetilde f$ as a section of $T^*U_x$ and denote its image
by $\operatorname{Im}(d\widetilde f)$.

\begin{proposition}
Let $X\subset U$ be a closed analytic subset and let $f:X\to\C$ be
holomorphic. A point $x\in X$ is an isolated critical point of $f$ if
there exist a neighborhood $U_x$ of $x$ and a holomorphic extension
$\widetilde f:U_x\to\C$ such that
$
\operatorname{Im}(d\widetilde f)\cap T_X^*U_x
=
\{(x,d\widetilde f(x))\}.
$
\end{proposition}

The next definition is taken from \cite[4.1]{BR}; compare with
\cite[p.~52]{GM} in the real case.

\begin{definition}\label{defn:Morse}
A holomorphic function $f:X\to\C$ is a \emph{(complex) Morse function}
if the following conditions hold:
\begin{enumerate}
\item For every critical point $x\in S_\alpha$ with
$\dim S_\alpha\geq 1$, the restriction $f|_{S_\alpha}$ has a
nondegenerate quadratic singularity at $x$.

\item For any holomorphic extension $\widetilde f$, any sequence
$x_n\in S_\beta$ converging to a critical point $x\in S_\alpha$, and
any limit plane $Q=\lim T_{x_n}S_\beta$, one has
$
d\widetilde f(x)(Q)\neq 0.
$
\end{enumerate}
\end{definition}

\begin{remark}
\begin{enumerate}
\item Condition (1) is equivalent to the transversality of the graph of
$d\widetilde f$ to the strata of the conormal
space. In particular, all critical points are isolated.

\item If $f$ is a complex Morse function, then $|f|$ is a real
stratified Morse function in the sense of \cite{GM}.

\item Condition (2) implies that all critical points of $|f|$ are
non-depraved in the sense of \cite{GM}.
\end{enumerate}
\end{remark}

\begin{definition}\label{defn}
A \emph{Morsification} of $f:X\to\C$ is a holomorphic family
$
F:X\times\D\to\C
$
such that $f_0=f$ and, for every $t\neq0$, the function $f_t$ is
Morse.
\end{definition}

 \begin{proposition}
If $f:(X,\0)\to(\C,0)$ has an isolated critical point, then it admits
a Morsification.
\end{proposition}

See \cite[Proposition 4.3]{BR}.

\begin{remark}
\begin{enumerate}
\item The Morsification may be chosen in the form
$
F(x,t)=f(x)+t\ell(x),
$
where $\ell$ is a general linear form (see the proof of \cite[Proposition 4.3]{BR}).

\item In particular, the restriction of a general linear form to a stratified
analytic set is Morse.

\item If $F$ is a Morsification, then the number of critical points
of $f_t$ near $x\in S_\alpha$ equals the Milnor number of
$f|_{S_\alpha}$ at $x$. This is a consequence of the conservation of the Milnor number (see e.g. \cite[Proposition E.2]{Mond-NB}).
\end{enumerate}
\end{remark}

\begin{example}
Let
$
X=\{(x,y,z)\in\C^3\mid xy=0\},
$
with strata
$
S_1=\{x=y=0\},\quad
S_2=\{x=0,\ y\neq0\},\quad
S_3=\{y=0,\ x\neq0\},
$
and let
$
f(x,y,z)=x^2+y^2+z^2.
$
Then $\0\in S_1$ is the only critical point, and $f|_{S_1}=z^2$ is
nondegenerate. However, if $x_n\in S_2$ tends to $\0$, then
$
\lim T_{x_n}S_2=\{0\}\times\C^2\neq T_{\0}S_1,
$
while $d\widetilde f(\0)=0$ for any extension $\widetilde f$. Thus,
condition (2) of Definition~\ref{defn:Morse} fails.
\end{example}


\section{Relative Polar Curves}

This section reviews standard results on polar curves relative to stratifications and linear forms. The properties needed in the sequel are summarized in Propositions \ref{good-forms-1} and \ref{good-forms-2}.


Let $X$ be an analytic subset of an open neighborhood $U$ of $\0$ in
$\C^N$, and let
$$
f:(X,\0)\to(\C,0)
$$
be the germ of a complex analytic function. We also denote by
$f:X\to\C$ a representative of this germ.

Let $\cS = (S_\a)_{\a \in A}$ be a Whitney stratification of a sufficiently small representative of $(X,\0)$. In particular, each stratum $S_\a$ and its closure $\overline{S}_\a$ are complex analytic spaces. We assume that $X$ is chosen small enough so that $\0\in\overline{S}_\a$ for all $\a\in A$, and that $A$ is finite. A stratification is said to be adapted to a complex analytic subspace if the latter is a union of strata.

We assume that the stratification is adapted to $f^{-1}(0)$. Then (see \cite{LM} or \cite[Theorem 6.7.1]{LNS}) there exists a non-empty Zariski open dense subset $\Omega_\a$ of the space of linear functions such that, for every $\ell\in\Omega_\a$ with $\ell(\0)=0$, the critical locus $C_\a$ of $(\ell,f)|_{S_\a\setminus f^{-1}(0)}$ is either empty or a smooth curve. Its closure $\Gamma_\a$ in $X$ is then either empty or a reduced curve. Moreover, the sets $\Omega_\a$ can be chosen so that the restriction $(\ell,f)|_{C_\a\setminus f^{-1}(0)}$ is finite for all strata.

Since the number of strata is finite, the intersection $\Omega:=\cap_{\a\in A}\Omega_\a$ is again a Zariski open dense subset. For any $\ell\in\Omega$, the union
$
\Gamma_\ell=\bigcup_{\a\in A}\Gamma_\a 
$
is either empty or a reduced curve.
Define (cf. \cite[p. 310]{L2}):

\begin{definition}  {\rm   For every $\ell\in\Omega_{\cS}:=\cap_{\a\in A}\Omega_\a$, the curve
$\Gamma_\ell=\Gamma_\ell(f,\cS,\0)$ defined above is called the
\emph{polar curve} of $f$ at $\0$ relative to $\ell$ and to the
stratification $\cS$.}
\end{definition}

  
 Let
$
\Delta_\ell=\Delta_\ell(f,\cS,\0):=\Phi(\Gamma_\ell)$, where $\Phi=(\ell,f).
$


\begin{definition}{\rm
If $\Delta_\ell(f,\cS,\0)$ is empty or a curve, it is called the
\emph{Cerf diagram} of $f$ at $\0$ relative to $\ell$ and to the
stratification $\cS$.}
\end{definition}

When no confusion is possible, we denote by $\Gamma_\ell$ the polar curve of $f$ at $\0$, and by $\Delta_\ell=\Delta_\ell(f,\cS,0)$ its Cerf diagram. If the stratification is fixed, we simply write $\Gamma_\ell=\Gamma_\ell(f,\0)$ and $\Delta_\ell=\Delta_\ell(f,\0)$.

By \cite[Theorem 6.7.1]{LNS}, for each $\a$ there exist $\varepsilon>0$ sufficiently small and a Zariski open dense subset $\Omega^1_\a$ of linear forms such that, for every $\ell\in\Omega^1_\a$, the map
$
\Phi=(\ell,f):X\cap \mathring{\mathbb B}_{\varepsilon'}\to\C^2,\quad 0<\varepsilon'\le \varepsilon,
$
is injective on $\Gamma_\a\cap \mathring{\mathbb B}_{\varepsilon'}$.

Moreover, the contribution of each stratum $S_\a$ to $\Gamma_\ell$ is either empty or a component of the polar curve (see \cite[Lemma 21]{LM} or \cite[Theorem 6.7.1]{LNS}). Since the index set $A$ is finite, one has:

\begin{proposition}\label{good-forms-1}
 There is a Zariski dense open set 
 $\Omega^1 = \Omega^1 (\cS, f) $ of linear forms such that for each $\ell  \in \Omega^1$ one has:
 \begin{enumerate}
 \item The critical locus of $\Phi$ is  a curve $\Gamma_\ell$ and its 
 restriction to each stratum that contains $\0$ in its closure is either empty or a curve. 
  \item The image $\Delta_\ell := \Phi(\Gamma_\ell)$ is either empty or a curve  in $\C^2$, called the Cerf diagram of $f$.
 \item The restriction of  $\Phi$ to $\Gamma_\ell$  is one-to-one. 
    \end{enumerate}
\end{proposition}

\begin{remark} {\rm
The above discussion extends to the case where the linear form $\ell$ is replaced by a holomorphic germ $g:(X,\0)\to(\C,0)$ with an isolated critical point at $\0$ relative to the given stratification. For each stratum $S_\a$, let $C_{g,f,\a}$ denote the critical set of $(g,f)|_{S_\a}$. It consists of points where either $f$ or $g$ is critical, or where their fibers are tangent in $S_\a$. The previous results remain valid under suitable genericity assumptions on these tangencies. This will be used later.}
\end{remark}

\subsection{Special choice of the linear form}\label{form}

We use stratified Morse theory to study the topology of the Milnor fiber.
This requires linear forms transverse to the limits of tangent hyperplanes.
For this purpose, we use the projectivized conormal space associated with
\ref{d. conormal}, which we denote by $C(X)$. More precisely,

$$C(X) \, = \overline{\{(x,H') \in X \times {\mathbb P}^{N-1} \, | \,   x \in  X_{\rm reg}   \,  \hbox{and \,}   T_x X \subset H'\}}\, . $$
where $X_{\rm reg}$ is the regular part of $X$. The projection
$\nu:C(X)\to X$ is given by $(x,H')\mapsto x$.
The fiber over $\0$ consists of the limits of hyperplanes containing
$T_xX$ as $x\in X_{\rm reg}$ tends to $\0$. These are called
\emph{tangent hyperplanes}.

By \cite{LT}, the space of limits of tangent hyperplanes to $X$ at $\0$ has dimension at most $N-2$, since it coincides with the space of hyperplanes tangent to the cones of the aureole of $X$ at $\0$ (see \cite[Definition 2.1.4, p.~562]{LT1}). 

For any stratum $S_\a$ with $\0\in\overline{S}_\a$, the set
$\Omega_{\overline{S}_\a,\0}$ of hyperplanes in $\C^N$ that do not
contain any limit of tangent spaces to $\overline{S}_\a$ at $\0$
is Zariski open and dense.
 Define
\[
\Omega^2=\bigcap_{S_\alpha\neq\{\0\}}\Omega_{\overline{S}_\alpha,\0},
\]
which is again a Zariski open dense subset of the space of linear forms on ${\mathbb C}^N$. We then have:

\begin{lemma}\label{transversality}
For every $\ell\in\Omega^2$, there exists an open neighborhood $V$ of
$\0$ in $X$ such that the restriction of $\ell$ to every stratum
$S_\alpha\cap V$ is a submersion, and the hyperplane $\ker\ell$
contains no limit of tangent spaces to these strata at $\0$.
\end{lemma}


Recall the notion of a stratified critical point relative to a Whitney stratification (Section~1.1). We will use the following observation:

\begin{lemma}
Let $X$ and $f$ be as before and let $\ell\in\Omega^2$ as in the
preceding subsection. Let $x$ be a point in
$X\cap\rond{\mathbb B}_\varepsilon\cap\Gamma_\alpha$. Then the
restriction of $\ell$ to the stratified space
$f^{-1}(f(x))\cap\rond{\mathbb B}_\varepsilon$
has an isolated stratified critical point at $x$.
\end{lemma}

\begin{definition} {\rm (See \ref {L0} below. Also,  {\it e.g.}, 
\cite[Theorem 6.11.3]{LNS})
We call the \emph{Milnor fiber} of $f$ at $\0$ the space
$f^{-1}(t)\cap\mathbb B_\varepsilon$, where
$0<|t|\ll\varepsilon\ll1$.}
\end{definition} 

\begin{lemma} \label{good forms}
For every $\ell\in\Omega^2$, the germ at $\0$ of
$\Phi=(\ell,f)$ satisfies the Thom condition.
\end{lemma}

See \cite[Theorem 2.6]{GLN} for a proof of Lemma~\ref{good forms}. 

Summarizing one has:

\begin{proposition}\label{good-forms-2}
There is a Zariski dense open set $\Omega^2$ of linear forms such that,
for each $\ell\in\Omega^2$, there exists an open neighborhood $V$ of
$\0$ in $X$ such that:
\begin{enumerate}
\item The restriction of $\ell$ to every stratum $S_\alpha\cap V$ is
a submersion, and the hyperplane $\ker\ell$ contains no limit of tangent
spaces to these strata at $\0$.

\item If $x\in V\cap\Gamma_\ell$, then the restriction of $\ell$ to the
stratified space $f^{-1}(f(x))\cap V$ has an isolated stratified
critical point at $x$.

\item The germ at $\0$ of the map $\Phi=(\ell,f)$ satisfies the Thom
condition.
\end{enumerate}
\end{proposition}

\begin{convention}\label{general forms}
{\rm The linear forms that we will consider from 
now on are those in $\Omega :=  \Omega^1 \cap \Omega^2$. 
This is a Zariski dense open subset of the space of all linear forms in $\C^N$, 
and each $\ell \in \Omega$ satisfies all properties in \ref{good-forms-1} and \ref{good-forms-2}. }
\end{convention}

\section{Preparation Lemmas}

Let $\cS=(S_\a)_{\a\in A}$ be a Whitney stratification of a sufficiently
small representative of $(X,\0)$, adapted to $f^{-1}(0)$, and such that
$\0\in\overline{S}_\a$ for all $\a\in A$. Choose a linear form $\ell\in\Omega:=\Omega^1\cap\Omega^2$ as above. Then the polar curve of $f$ at $\0$ relative to $\ell$ and $\cS$ is defined.

The induced stratification of $\ell^{-1}(0)\cap f^{-1}(0)$ is given by $(S_\a\cap \ell^{-1}(0)\cap f^{-1}(0))$.

Let $\varepsilon>0$ be sufficiently small so that the sphere 
$\mathbb S_{\varepsilon}:=\partial \mathbb B_{\varepsilon}$ is a Milnor sphere for $X$, $f^{-1}(0)$, and $f^{-1}(0)\cap \ell^{-1}(0)$. Thus, for any $0<\varepsilon'\le \varepsilon$, the sphere $\mathbb S_{\varepsilon'}$ intersects transversely the strata of these spaces. This follows from the local conic structure of analytic sets (see \cite[Lemma 3.2]{BV}, cf.\ Lemma \ref{classic} below).

Under these assumptions, we have:

\begin{lemma}\label{Thom} The restriction $\ell|_{f^{-1}(0)}$ has an isolated stratified critical
point at $\0$ with respect to the stratification $(S_\alpha)$, and there
exists $\eta>0$, with $1\gg\varepsilon\gg\eta$,  such that, for 
any $(u,t)\in \rond {\D}_\eta\subset {\mathbb C}^2,$ the fiber $\ell^{-1}(u)\cap f^{-1}(t)$ is transverse 
to ${\mathbb S}_\varepsilon$  and has only isolated stratified critical points of $\ell_{|f^{-1}(t)}$
with respect to the stratification $(S_\alpha\cap f^{-1}(t))$.
\end{lemma}

The proof is straightforward. 

We may assume that $\eta$ is chosen so that
$\Delta_\ell\cap\{t=0\}\cap\rond{\D}_\eta=\{0\}.$
If $\Delta_\alpha$ is empty for every $\alpha\in A$, then, using
rugose vector fields, one can prove that
$\;X\cap\mathbb B_\varepsilon$ is homeomorphic to
$
X\cap\ell^{-1}(0)\cap\mathbb B_\varepsilon\times\mathbb D_\eta \,
$
(see Lemma~\ref{level a} below).

\begin{definition}\label{delta}
We denote by $\delta_\a=(\Delta_\a\cdot\{t=0\})_0$ the intersection
number of the curves $\Delta_\alpha$ and $\{t=0\}$ at the origin.
One can choose $|t_0|$ small enough so that the line $\{t=t_0\}$
intersects the component $\Delta_\alpha$ at distinct points contained
in the open ball $\rond{\D}_\eta$. Then $\delta_\a$ is the number of
such points.
\end{definition}

Thus $\delta_\alpha$ is the number of stratified critical points of
$\ell$ on $S_\alpha\cap f^{-1}(t_0)$.

The following lemma is well-known and we omit the proof:

\begin{lemma}\label{classic}
For $\varepsilon>0$ small enough, the neighborhood
$\mathbb B_\e\cap X$ of $\0$ in $X$ is good in the sense of Prill
\cite{Prill}: it is stratified homeomorphic to the cone over its link
$L_X:=\mathbb S_\e\cap X$, and, after removing the vertex $\0$, it is
stratified homeomorphic to the cylinder $L_X\times[0,1)$. Also, for
$0<|t_0|\ll\varepsilon$, $f^{-1}(t_0)\cap\mathbb B_\e$ is homeomorphic
to the Milnor fiber of $f$ at $\0$.
\end{lemma}

\begin{notation}\label{notation} We let $a,b \in \R$,  $0<a<b$ , 
be such that all the points of $\Delta_\alpha\cap \{t=t_0\}\cap {\D}_\eta$ lie in the region
$a<|\ell|<b$.
\end{notation}

The proof of the next lemma is a straightforward application of the
fibration theorem in \cite{L0} (see \ref{L0} below) and of rugose vector
fields adapted to a Whitney-strong (or $w$-regular) stratification
(see \cite{V}). These tools will also be used in the proof of
Lemma~\ref{Milnor}.

In \cite[p.~395]{T}, Teissier introduced $w$-regularity under the name \emph{Condition (a) stricte avec exposant 1}, and showed that, in the complex analytic setting, it is equivalent to the Whitney conditions (see \cite[Theorem 1.2, p.~455]{T}). Verdier used such $w$-regular stratifications in \cite{V} to define and construct \emph{rugose vector fields}, that is, stratified continuous vector fields that are smooth on each stratum and integrable.

\begin{lemma}\label{level a}
Let $a$ be as in \ref{notation} and let $B_a$ be the ball in $\C$
centered at $0$ with radius $a$. Then the subset
$f^{-1}(t_0)\cap\mathbb B_\varepsilon\cap\ell^{-1}(B_a)$ of the
Milnor fiber $F_f:=f^{-1}(t_0)\cap\mathbb B_\varepsilon$ is
homeomorphic to $F_{\ell,f}\times B_a$ by a stratum-preserving
homeomorphism that is smooth on each stratum.
\end{lemma}

\begin{figure}[H]
\begin{center}
\includegraphics[width=7.5cm]{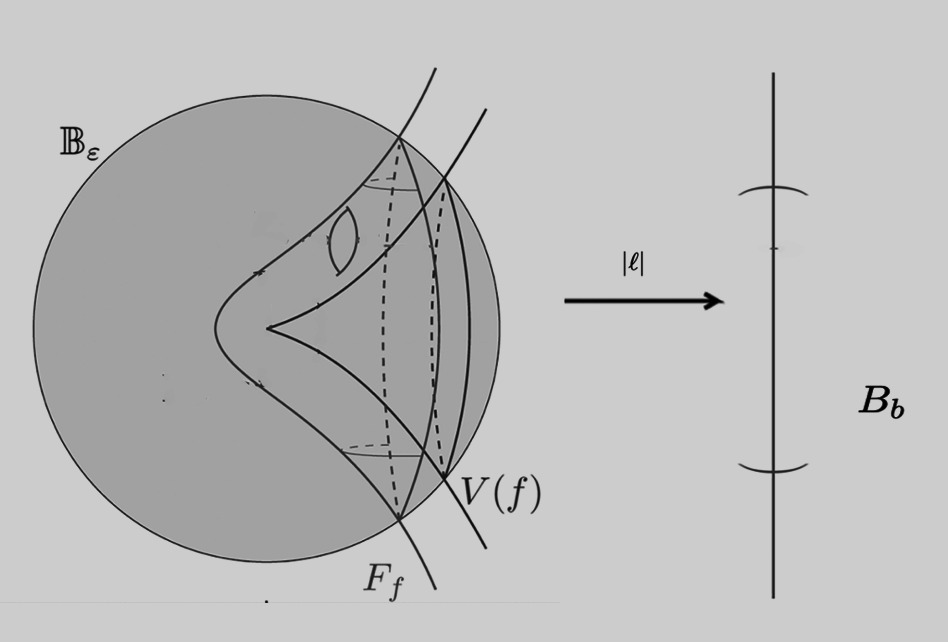}\label{Fig 1}
\caption{}\label{fig:LNS}
\end{center}
\end{figure}

Next we have (see Figure \ref{fig:LNS}):

\begin{lemma}\label{Milnor}
Let $b$ be as in \ref{notation} and let $B_b$ be the ball in $\C$
centered at $0$ with radius $b$. Then the subset
$f^{-1}(t_0)\cap\mathbb B_\varepsilon\cap\ell^{-1}(B_b)$ of the
Milnor fiber $F_f:=f^{-1}(t_0)\cap\mathbb B_\varepsilon$ is
homeomorphic to the whole fiber $F_f$ by a stratum-preserving
homeomorphism that is smooth on each stratum.
\end{lemma}

\begin{proof}
To prove Lemma~\ref{Milnor}, we proceed as in Milnor
\cite[Chap.~4]{M}. Since $X$ may have non-isolated singularities, the
Milnor fiber itself can be singular, so we construct a rugose vector
field whose integration gives the required homeomorphism.

Consider the following open subset (see Figure~\ref{fig:bolas}):
$$
A(\varepsilon',b')
:=
f^{-1}(t_0)\cap\rond{\mathbb B}_{\varepsilon'}
-\ell^{-1}(B_{b'}),
$$
where $0<\varepsilon'-\varepsilon\ll1$ and $0<b-b'\ll1$.

\begin{figure}[H]
\begin{center}
\includegraphics[width=10cm]{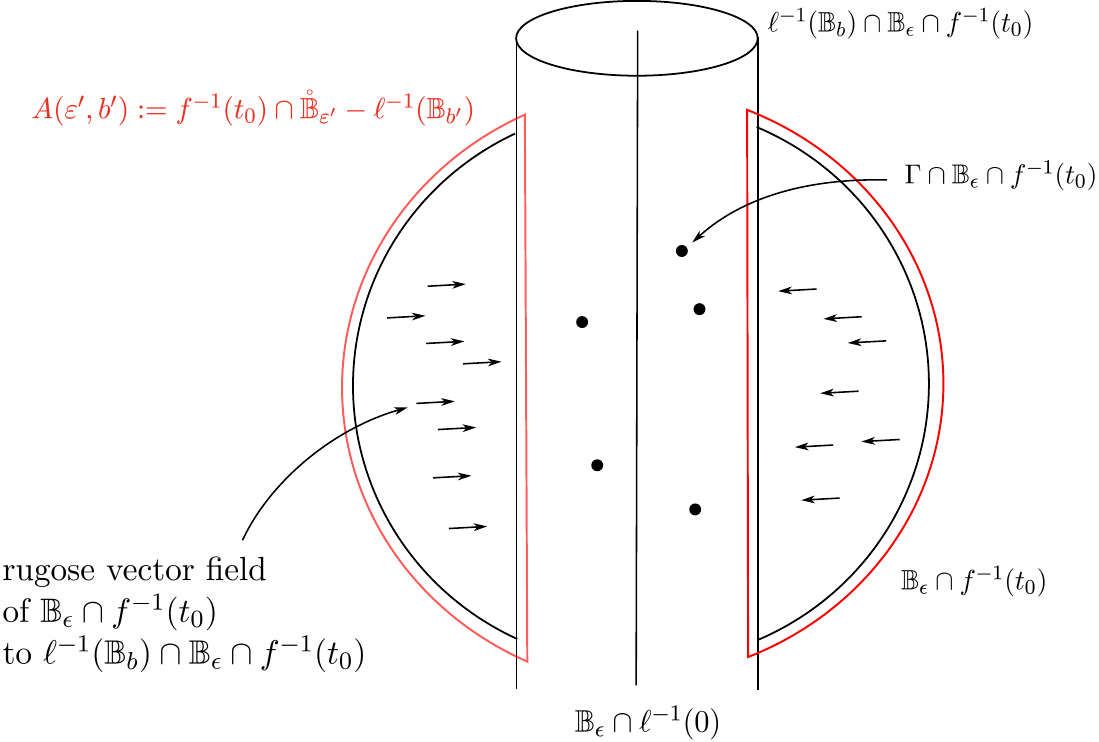}
\caption{}\label{fig:bolas}
\end{center}
\end{figure}

The Whitney stratification $(S_\alpha)$ of $X$, adapted to the fiber
$f^{-1}(0)$, induces a Whitney stratification of
$A(\varepsilon',b')$ (see Theorem~2.6 of \cite{GLN}). By the choice
of $b$, for suitable $\varepsilon'$ and $b'$, the space
$A(\varepsilon',b')$ does not intersect the polar curve $\Gamma$.
Its strata are
$$
T_\alpha:=S_\alpha\cap A(\varepsilon',b'),
$$
and have dimension at least $1$. Fix such $\varepsilon'$ and $b'$.
There are only finitely many indices $\alpha$ for which
$T_\alpha\neq\emptyset$.

Let $grad_\alpha\|x\|$ and $grad_\alpha|\ell|(x)$ be the gradients of
the restrictions of $\|x\|$ and $|\ell|$ to $T_\alpha$. Let
$x\in T_\alpha\cap\mathbb B_\varepsilon$. Then, by the Curve Selection
Lemma, as in \cite{M}, either:
\begin{enumerate}
\item the vectors $grad_\alpha\|x\|$ and
$grad_\alpha|\ell|(x)$ are linearly independent; or
\item the vectors $grad_\alpha\|x\|$ and
$grad_\alpha|\ell|(x)$ are collinear, and
$$
grad_\alpha\|x\|
=
\lambda_x grad_\alpha|\ell|(x)
$$
for some $\lambda_x\in\R$.
\end{enumerate}

In the first case, there is an open neighborhood $U_x$ of $x$ in
$A(\varepsilon',b')$ and a rugose stratified vector field $v_x$
(see \cite[S4]{V}) relative to the stratification $(T_\alpha)$ such
that, at any point $y\in T_{\alpha_y}\cap U_x$, for some $\alpha_y$, the vector $v_x(y)$
is tangent to $T_{\alpha_y}$ and satisfies
$$
Re\langle v_x(y),grad_{\alpha_y}\|y\|\rangle<0
$$
and
$$
Re\langle v_x(y),grad_{\alpha_y}|\ell|(y)\rangle<0,
$$
where $Re\langle\, ,\,\rangle$ is the real part of the Hermitian
product in $\C^N$.

In the second case, if
$$
grad_\alpha\|x\|
=
\lambda_x grad_\alpha|\ell|(x),
\qquad \lambda_x\in\R,
$$
then $\lambda_x\geq0$ when $x$ belongs to a sufficiently small
neighborhood $U_0$ of $\0$ in $X$ (see
\cite[Corollary~3.4]{M}). We may assume that
$\mathbb B_\varepsilon\subset U_0$ and that
$\lambda_x\neq0$ on $U_0\cap A(\varepsilon',b')$.

We may choose $\varepsilon'$ and $b'$ such that
$A(\varepsilon',b')$ contains
$$
A:=F_f\setminus\ell^{-1}(\rond{B}_b)
$$
and is contained in $U_0$.

Thus, for every $x\in A(\varepsilon',b')$, there is an open
neighborhood $U_x\subset A(\varepsilon',b')$ and a nowhere vanishing
rugose vector field $v_x$ on $U_x$ such that, for every
$y\in T_{\alpha_y}\cap U_x$,
$$
Re\langle v_x(y),grad_{\alpha_y}\|y\|\rangle<0
$$
and
$$
Re\langle v_x(y),grad_{\alpha_y}|\ell|(y)\rangle<0.
$$

After multiplying the local vector fields $v_x$ by suitable positive
constants and shrinking the neighborhoods $U_x$ if necessary, we may
assume that
$$
Re\langle v_x(y),grad_{\alpha_y}|\ell|(y)\rangle\leq -1
$$
on $U_x$.

Let $(\varphi_x)$ be a rugose partition of unity subordinate to the
covering $\bigcup_x U_x$ of $A(\varepsilon',b')$. Then
$$
v:=\sum_x\varphi_xv_x
$$
is a rugose vector field on $A(\varepsilon',b')$ satisfying
$$
Re\langle v(y),grad_{\alpha_y}\|y\|\rangle<0
$$
and
$$
Re\langle v(y),grad_{\alpha_y}|\ell|(y)\rangle\leq-1.
$$

Choose a nonnegative smooth function $\Psi$ on
$A(\varepsilon',b')$, equal to $1$ on a neighborhood of $A$ and with
compact support in $A(\varepsilon',b')$. Then $\Psi v$ is a rugose
vector field with compact support in $A(\varepsilon',b')$ and agrees
with $v$ on a neighborhood of $A$.

Since all the critical points of $|\ell|$ on $F_f$ lie in the region
$a<|\ell|<b$, we may choose $c<b$ such that
$$
\max\{|\ell(p)|: p\text{ is a critical point of }|\ell|_{|F_f}\}
<c<b.
$$
Thus $|\ell|$ has no critical points on
$$
F_f\cap\{|\ell|\geq c\}.
$$
The inequalities above imply that the flow remains in the interior of
$\mathbb B_\varepsilon$ and crosses the regular levels of $|\ell|$ in
finite time. Hence the flow gives a stratum-preserving product
structure on
$$
F_f\cap\{c\leq|\ell|\leq M\},
$$
where
$
M:=\max_{F_f}|\ell|.
$
Choose an increasing homeomorphism
$$
h:[c,M]\longrightarrow[c,b]
$$
such that $h(c)=c$ and $h(r)\leq r$ for every $r\in[c,M]$.

Moving each point along its flow line from the level
$|\ell|=r$ to the level $|\ell|=h(r)$, and leaving
$
F_f\cap\{|\ell|\leq c\}
$
fixed, gives a stratum-preserving homeomorphism
$$
F_f\longrightarrow F_f\cap\ell^{-1}(B_b)
$$
which is smooth on each stratum. This proves Lemma~\ref{Milnor}.

\end{proof}





\section{Stratified Morse Theory}\label{s: theorem linear case}
Choose a linear form $\ell\in\Omega^1\cap\Omega^2$ as in
\ref{general forms}. We assume that $\varepsilon>0$ is sufficiently
small so that the restriction to the polar curve $\Gamma$ of the map
$\Phi:=(\ell,f)$ introduced in Section~\ref{form} is one-to-one onto
the Cerf diagram $\Delta$.

Consider the real valued function $|\ell|$ and  its restriction  to the Milnor fiber  
$F_f := f^{-1}(t_0)\cap {\mathbb B}_\varepsilon$.

 \begin{lemma}\label{critical}
The stratified critical points of the restriction
$|\ell|_{|F_f\cap|\ell|^{-1}((a,b))}$, a real valued function,
are the stratified critical points of the complex valued function
$\ell_{|F_f\cap|\ell|^{-1}((a,b))}$.
\end{lemma}

The next result implies that with our choice of the linear form $\ell$, the critical points of 
$|\ell|_{|F \cap \ell^{-1}((a,b))}$ are Morse singularities: 

 \begin{lemma}\label{ordinary}
A critical point $x$ of the restriction of $\ell$ to a stratum $S_\alpha$ is an ordinary quadratic 
point of the restriction $|\ell|_{|S_\alpha}$,  
if $\ell(x)\neq 0$.
\end{lemma}

\begin{proof}

In the proof of Theorem 2.6 of \cite{GLN}, it is shown that the Whitney
stratification of $X$, adapted to $f^{-1}(0)$, can be refined so as to
be adapted to the polar curve $\Gamma$ of $f$ relative to $\ell$.
Intersecting this stratification with
$F_f\cap|\ell|^{-1}((a,b))$ yields a Whitney stratification of
$F_f\cap|\ell|^{-1}((a,b))$ whose strata are
$$
F_f\cap S_\alpha-\bigcup_{\alpha,j}\{z_{\alpha,j}\}
$$
and the points $\{z_{\alpha,j}\}$, for $\alpha$ and
$1\leq j\leq j_\alpha$.

The stratified critical points of the restriction
$$
|\ell|_{|F_f\cap|\ell|^{-1}((a,b))}
$$
are the points of $F_f\cap\ell^{-1}((a,b))$ where $|\ell|$ is not
transverse to the strata of
$F_f\cap S_\alpha\cap\ell^{-1}((a,b))$ induced by $(S_\alpha)$.
The only $0$-dimensional strata are the points $\{z_{\alpha,j}\}$,
while all other strata have dimension at least $1$, and $|\ell|$ is
transverse to them. Therefore, the only stratified critical points of
$$
|\ell|_{|F_f\cap|\ell|^{-1}((a,b))}
$$
are the points $\{z_{\alpha,j}\}$, namely the intersections of the
polar curve $\Gamma$ with $F_f$.

The critical set of
$$
\Phi_{|S_\alpha\cap\mathring{\mathbb B}_\varepsilon}
=
(\ell_{|S_\alpha\cap\mathring{\mathbb B}_\varepsilon},
f_{|S_\alpha\cap\mathring{\mathbb B}_\varepsilon})
$$
on the stratum $S_\alpha$ is the reduced curve
$$
\Gamma_\alpha\cap\mathring{\mathbb B}_\varepsilon-\{\0\}.
$$
Hence $\Phi$ has ordinary quadratic singularities at points of
$\Gamma_\alpha\cap\mathring{\mathbb B}_\varepsilon$. 
It follows that the restriction of $\ell$ to
$S_\alpha\cap F_f$ has ordinary critical points at
$\{z_{\alpha,j}\}$, for $\alpha$ and $1\leq j\leq j_\alpha$
(see \cite[Theorem 6.7.6]{LNS}).

\end{proof}

Now let $a$ and $b$ be as in \ref{notation}, so that the open interval
$(a,b)$ contains all $|\ell(z_{\alpha,j})|$, where the
$z_{\alpha,j}$ are the intersection points of the polar curve
$\Gamma$ with the Milnor fiber $F_f=f^{-1}(t_0)\cap\mathbb B_\varepsilon$.
Equivalently, the points
$
(\ell(z_{\alpha,j}),t_0)
$
are the intersection points of the line $t=t_0$ with the Cerf diagram
$\Delta$ in $\mathbb D_\eta\subset\mathbb C^2$.
 Then we have:

\begin{lemma}\label{Morse}
With $\ell$, $\eta$, and $t_0$ as above, the function
$
|\ell|_{|\{f=t_0\}\cap\mathbb B_\varepsilon}
$
is a stratified real Morse function on
$
f^{-1}(t_0)\cap\mathbb B_\varepsilon\cap |\ell|^{-1}((a',b')),
$
for $0<a'\leq a<b\leq b'\leq b+\delta$, where $\delta>0$ is sufficiently small.
\end{lemma}
Before proving the lemma, recall that a map $X\to\mathbb R$ is Morse if it is proper, has only isolated nondegenerate critical points on each stratum, and if, at each critical point $x\in S_\a$, its differential does not vanish on any limit of tangent spaces to strata $S_\beta\ne S_\a$ (see, for instance, \cite{GM}).

\begin{proof}
First observe that the restriction
$\,
|\ell|_{|f^{-1}(t_0)\cap\{a<|\ell|<b\}\cap\mathbb B_\varepsilon}
$
is proper.

Let $\partial F_f:=F_f\cap\mathbb S_\varepsilon$. The strata of $F_f-\partial F_f$ are the intersections of the strata $S_\alpha$ of $X$ with $F_f-\partial F_f$, while the strata of $\partial F_f$ are the intersections of the strata $S_\alpha$ with $\partial F_f$.

By the choice of $\varepsilon>0$ and $\eta$ in the previous section (Lemma \ref{Thom}), the form $\ell$ is transverse to the strata $S_\alpha\cap\partial F_f$ of $\partial\mathbb B_\varepsilon\cap f^{-1}(t)$ for every $t$ such that $(u,t)\in\mathbb D_\eta$. Hence the only stratified critical points of the restriction of $|\ell|$ to $F_f\cap\{a<|\ell|<b\}$ lie in $F_f-\partial F_f$.

The stratified critical points of the restriction of $|\ell|$ to
$
(F_f-\partial F_f)\cap\{a<|\ell|<b\}
$
are precisely the points where $F_f$ meets the polar curve
$$
\Gamma=\bigcup_\alpha \Gamma_\alpha
$$
of $f$ at $\0$. Since each $\Gamma_\alpha$ is nonsingular outside $\0$ (see \cite[Lemma 21 (1)]{LM} or \cite[Theorem 6.7.6]{LNS}), the restriction of $\ell$ to $S_\alpha\cap\mathring{\mathbb B}_\varepsilon$ is ordinary quadratic at every point where $\ell_{|S_\alpha}$ is critical. By Lemma \ref{ordinary}, the restriction of $|\ell|$ to $S_\alpha$ is therefore real ordinary quadratic at every point of $S_\alpha\cap\Gamma_\alpha$.

The points $z_{\alpha,j}$, $j=1,\ldots,j_\alpha$, in $F_f\cap\Gamma_\alpha$ are the stratified critical points of the restriction of $|\ell|$ to $F_f$, and satisfy
$
a<|\ell|(z_{\alpha,j})<b.
$

It remains to verify the third condition in the definition of a stratified Morse function (see \cite[Definition 5.5.2]{G}). Since $\mathbb B_\varepsilon\subset V$ (see Lemma \ref{transversality}), the linear form $\ell$ is transverse in $V$ to the limits of tangent spaces to $\overline{S}_\a$ at $z_{\a,j}$. Hence, at each $z_{\a,j}$, the tangent space of $|\ell|$ restricted to the corresponding stratum is transverse to all limits of tangent spaces to strata whose closures contain $z_{\a,j}$.

By \S 5.5.2 of \cite{G}, it follows that the restriction of $|\ell|$ to $F_f$ is a non-depraved stratified Morse function on
$
f^{-1}(t_0)\cap\mathbb B_\varepsilon\cap |\ell|^{-1}((a,b)).
$

The choice of $a'$ and $b'$ follows from the fact that
$$
0<a'\leq a<b\leq b'\leq b+\delta,
$$
and that all the points $z_{\alpha,j}$ lie in
$
\{a'\leq|\ell|\leq b'\}.
$
This proves Lemma~\ref{Morse}.
\end{proof}

We now apply stratified Morse theory to the function $|\ell|_{|F_f}$,
using it to pass from
$
F_f\cap\{|\ell|\leq a\}
$
to
$
F_f\cap\{|\ell|\leq b\}.
$

By Thom's First Isotopy Lemma \cite[\S 11]{Mather} (see also
\cite[Theorem 5.5.1]{G}), if $a_1>a$ and there are no critical points
of $|\ell|_{|F_f}$ with $a<|\ell|\leq a_1$, then
$
F_f\cap\{|\ell|\leq a\}
$
and
$
F_f\cap\{|\ell|\leq a_1\}
$
are homeomorphic by a stratum-preserving homeomorphism that is smooth
on each stratum.

The next step is to describe what happens when one crosses a critical
level of $|\ell|$ on $F_f$; we follow \cite{GM} (see
\cite[Theorem 5.5.3]{G}). Recall that the critical points of $|\ell|$
on $F_f$ are precisely the points where the polar curve $\Gamma$ meets
$F_f$. Thus:

\begin{lemma}\label{l:}
The space
$
f^{-1}(t_0)\cap\mathbb B_\varepsilon\cap\ell^{-1}(B_b)
$
is homeomorphic to the space obtained from
$
f^{-1}(t_0)\cap\mathbb B_\varepsilon\cap\ell^{-1}(B_a)
$
by attaching $\sum_\alpha j_\alpha$ Morse data, one for each critical
point $z_{\alpha,j}$ of
$
|\ell|_{|F_f\cap\{a<|\ell|<b\}}.
$
\end{lemma}

We now describe the corresponding Morse data. By definition, it is the product of the tangential and normal Morse data.

Since the critical points of $\ell$ are ordinary quadratic
singularities, let
$
d_\alpha:=\dim_{\C}(S_\alpha\cap F_f).
$
At each critical point $z_{\alpha,j}\in S_\alpha\cap F_f$, the
restriction of $|\ell|$ has Morse index $d_\alpha$. Hence the
tangential Morse data is
$
(D^{d_\alpha},\partial D^{d_\alpha})\times D^{d_\alpha},
$
where $D^{d_\alpha}$ denotes a closed disc of dimension $d_\alpha$
(see the remark following Proposition 3.6.2, p.~65 of \cite{GM}).

The normal Morse data is independent of the function and depends only
on the space $X$ and the stratum $S_\alpha$ containing
$z_{\alpha,j}$. Its total space is homeomorphic to a normal slice at
$z_{\alpha,j}$, hence to the cone over the link $L_\alpha$ of the
stratum $S_\alpha$ at $z_{\alpha,j}$.


To describe how this space is attached, one needs the \emph{half
links}, the real analogues of the complex link. There are an upper
half link $(L_\alpha^+,\partial L_\alpha^+)$ and a lower half link
$(L_\alpha^-,\partial L_\alpha^-)$; see \cite[3.9.1]{GM}. These are
stratified spaces whose union along their common boundary is
homeomorphic to the link:
$
L_\alpha=L_\alpha^+\cup_\partial L_\alpha^-.
$

By Theorem 3.11.1 of \cite{GM}, the normal Morse data at each critical
point is
$$
(\operatorname{Cone}(L_\alpha^+\cup_\partial L_\alpha^-),L_\alpha^-).
$$

Furthermore, in the complex analytic setting, the normal Morse data
admits a simpler homotopy description (see
\cite[Corollary 1, p.~166]{GM}): it has the homotopy type of the pair
$
(\operatorname{Cone}(\mathcal L_\alpha),\mathcal L_\alpha),
$
where $\mathcal L_\alpha$ is the complex link of the stratum.

Let us summarize. Given a Whitney stratified complex analytic map-germ
$$
f:(X,\0)\to(\C,0), \qquad X\subset\C^N,
$$
and its Milnor fiber
$
F_f=f^{-1}(t_0)\cap\mathbb B_\varepsilon,
$
with the induced Whitney stratification, there exists a Zariski dense
open set $\Omega$ of linear forms on $\C^N$ such that, for every
$\ell\in\Omega$, if we set $\Phi:=(\ell,f)$, then:

\begin{itemize}
\item The polar set of $(\ell,f)$ is a curve $\Gamma$ meeting
$\Phi^{-1}(0,0)$ only at the origin; $\Phi$ satisfies the Thom
condition, and its restriction to $\Gamma$ is one-to-one onto its
image.

\item The function $|\ell|$ restricted to
$F_f\cap|\ell|^{-1}((a',b'))$ is a nondepraved stratified Morse
function. Its critical points are precisely the intersection points
of $F_f$ with the polar curve $\Gamma$.

\item At each critical point $z_{\alpha,j}\in S_\alpha\cap F_f$,
$\ell$ has an ordinary quadratic singularity. If
$
d_\alpha:=\dim_{\C}(S_\alpha\cap F_f),
$
then the restriction of $|\ell|$ to $S_\alpha\cap F_f$ has Morse
index $d_\alpha$. The tangential and normal Morse data are,
respectively,
$$
(D^{d_\alpha},\partial D^{d_\alpha})\times D^{d_\alpha}
$$
and
$$
(\operatorname{Cone}(L_\alpha^+\cup_\partial L_\alpha^-),
L_\alpha^-),
$$
where $L_\alpha^\pm$ are the half links.

\item In the complex analytic setting, the normal Morse data is
homotopy equivalent to
$
(\operatorname{Cone}(\mathcal L_\alpha),\mathcal L_\alpha),
$
where $\mathcal L_\alpha$ is the complex link.

\item The level $\{|\ell|\leq a\}\cap F_f$ is homeomorphic to
$F_{\ell,f}\times\mathbb D^2$, where $F_{\ell,f}$ is the Milnor fiber
of $(\ell,f)$.

\item The level $\{|\ell|\leq b\}\cap F_f$ is homeomorphic to the
whole Milnor fiber $F_f$.

\end{itemize}

We get:

\begin{theorem}\label{main}
Let $F_f$ be the Milnor fiber of $f$ at $\0$, let $\ell$ be a
sufficiently general linear form as above, and let $F_{\ell,f}$ be the
Milnor fiber of $\Phi:=(\ell,f)$. For each stratum $S_\alpha$ of $X$,
set
$
d_\alpha:=\dim_{\C}(S_\alpha\cap F_f),
$
and recall that
$
\delta_\alpha=(\Delta_\alpha\cdot\{t=0\})_0.
$
Then $\delta_\alpha$ is the number of critical points of the
restriction of $|\ell|$ to $S_\alpha\cap F_f$.
The Milnor fiber $F_f$ is homeomorphic to
$
F_{\ell,f}\times\mathbb D^2
$
with, for each stratum $S_\alpha$, $\delta_\alpha$ copies of the Morse
data
$$
(D^{d_\alpha},\partial D^{d_\alpha})\times D^{d_\alpha}\times
(\operatorname{Cone}(L_\alpha^+\cup_\partial L_\alpha^-),
L_\alpha^-)
$$
attached.
Moreover, replacing each Morse datum by
$$
(D^{d_\alpha},\partial D^{d_\alpha})\times D^{d_\alpha}\times
(\operatorname{Cone}(\mathcal L_\alpha),\mathcal L_\alpha),
$$
where $\mathcal L_\alpha$ is the complex link of the stratum, one
recovers $F_f$ up to homotopy.
\end{theorem}

\begin{proof}
By Lemma \ref{Morse} and stratified Morse theory, the space
$
f^{-1}(t_0)\cap\mathbb B_\varepsilon\cap \ell^{-1}(B_b)
$
is obtained from
$
f^{-1}(t_0)\cap\mathbb B_\varepsilon\cap \ell^{-1}(B_a)
$
by attaching Morse data at each stratified critical point $x$ of
$|\ell|_{|F_f}$ such that
$
a<|\ell|(x)<b.
$

For each stratum $S_\alpha$, there are precisely $\delta_\alpha$
such critical points. At each critical point
$x\in S_\alpha\cap F_f$, the Morse data is
$$
(D^{d_\alpha},\partial D^{d_\alpha})\times D^{d_\alpha}\times
(\operatorname{Cone}(L_\alpha^+\cup_\partial L_\alpha^-),
L_\alpha^-).
$$
Hence
$
f^{-1}(t_0)\cap\mathbb B_\varepsilon\cap \ell^{-1}(B_b)
$
is homeomorphic to the space obtained from
$
f^{-1}(t_0)\cap\mathbb B_\varepsilon\cap \ell^{-1}(B_a)
$
by attaching, for each stratum $S_\alpha$, $\delta_\alpha$ copies of
the above Morse data.

By Lemma \ref{level a},
$
f^{-1}(t_0)\cap\mathbb B_\varepsilon\cap \ell^{-1}(B_a)
$
is homeomorphic to
$
F_{\ell,f}\times\mathbb D^2,
$
while, by Lemma \ref{Milnor},
$
f^{-1}(t_0)\cap\mathbb B_\varepsilon\cap \ell^{-1}(B_b)
$
is homeomorphic to $F_f$. This proves the first assertion.

Finally, in the complex analytic setting, the normal Morse data has
the homotopy type of
$
(\operatorname{Cone}(\mathcal L_\alpha),\mathcal L_\alpha).
$
Therefore, replacing the normal Morse data in each attachment by this
pair gives the corresponding description of $F_f$ up to homotopy.
This proves the second assertion.
\end{proof}


\section{Milnor-L\^{e} fibrations and relative Morsifications}\label{s: theorem general case}

Let $(X,\0)$ be a germ as above, where $X$ is a closed analytic subset of an open neighborhood $U\subset \C^N$ of $\0$, endowed with a Whitney stratification $(S_\a)$. Fix a Milnor ball $\B_\varepsilon$ for $X$ at $\0$, that is, for every $0<\varepsilon'\le \varepsilon$, the sphere $\s_{\varepsilon'}$ is transverse to all strata $S_\a$. In particular, this induces a Whitney stratification on $\B_\varepsilon\cap X$.

\subsection{The Thom  condition}
Suppose that $f=(f_1,\dots,f_k)\colon (X,\0)\to(\C^k,0)$ is a holomorphic map germ, and fix a representative $f\colon X\to V$, where $V\subset\C^k$ is 
an open neighborhood of the origin. A  {\it stratification} of $f$ consists of Whitney stratifications 
$(S_\alpha)$ of $X$ and $(T_\beta)$ of $V$ such that, for 
each $S_\alpha$, there exists $T_\beta$ with $f(S_\alpha)\subset T_\beta$ and $f\colon S_\alpha \to T_\beta$ a submersion. 

It is known that, after a suitable choice of strata, every map between real or complex analytic spaces admits such a stratification (see \cite{Durfee} for the semi-algebraic case; the analytic case is analogous).

The \emph{discriminant} of $f$ with respect to the stratification is defined as
$
\widehat \Delta:=\bigcup_{\dim T_\beta<k} \overline{T_\beta} \;. 
$
By construction, $\widehat\Delta$ is a closed analytic subset of $V$ of dimension $<k$. Moreover,
$$
f\colon X\cap f^{-1}(V\setminus\widehat\Delta)\longrightarrow V\setminus\widehat\Delta
$$
is a stratified submersion, and for every $u\in V\setminus\widehat\Delta$, the fibre $f^{-1}(u)$ is a stratified subspace of codimension $k$ in $X$ with the induced stratification. If $V$ is connected, then so is $V\setminus\widehat\Delta$. However, this map is generally not proper, so the Thom--Mather First Isotopy Lemma does not directly imply local $C^0$-triviality.

\begin{definition}\label{def:fib}{\rm
Let $f=(f_1,\dots,f_k)\colon (X,\0)\to(\C^k,0)$ be a holomorphic map germ. We say that $f$ admits a \emph{Milnor-L\^{e} fibration} if there exists a stratification of a representative $f\colon X\to V$ such that, for any $\epsilon,\delta$ with $0<\delta\ll\epsilon\ll1$, the restriction
\begin{equation}\label{ML}
f\colon \B_\epsilon\cap X\cap f^{-1}(\mathring{\D}_\delta\setminus\widehat\Delta)\longrightarrow \mathring{\D}_\delta\setminus\widehat\Delta
\end{equation}
is a locally trivial stratified fibration with the induced stratification.
}\end{definition}


\begin{lemma} \label{l: Milnor-Le fibration}
Assume that $f\colon X\to V$ admits a stratification which satisfies the Thom condition. Then $f$ admits a Milnor-L\^{e} fibration.
\end{lemma}

\begin{proof}  By the Thom--Mather First Isotopy Lemma \cite[\S 11]{Mather}, a proper stratified map is a locally trivial stratified fibration if and only if it is a stratified submersion.

For any $\epsilon,\delta>0$ and any stratum $S_\alpha\subset X$, the restriction of $f$ to
$
S_\alpha\cap\mathring{\B}_\epsilon\cap f^{-1}(\mathring{\D}_\delta\setminus\widehat\Delta)
$
is already a submersion. Thus it remains to choose $\epsilon,\delta>0$ so that $f$ is also a submersion on strata of the form
$
S_\alpha\cap\s_\epsilon\cap f^{-1}(\mathring{\D}_\delta\setminus\widehat\Delta).
$
Equivalently, the fibres of the map \eqref{ML} in Definition \ref{def:fib} must be stratified transverse to the sphere $\s_\epsilon$.

Choose first $\epsilon>0$ so that $\s_\epsilon$ is transverse to every stratum $S_\alpha\subset f^{-1}(0)$. We claim that there exists $\delta>0$ such that, for every
$
x\in \s_\epsilon\cap f^{-1}(\mathring{\B}_\delta\setminus\widehat\Delta),
$
the sphere $\s_\epsilon$ is transverse to
$
f^{-1}(f(x))\cap S_\beta,
$
where $S_\beta$ is the stratum containing $x$. In particular, $f$ is then a submersion on
$
S_\beta\cap\s_\epsilon\cap f^{-1}(\mathring{\D}_\delta\setminus\widehat\Delta).
$

If not, there would exist a sequence $\{x_n\}$ in
$
\s_\epsilon\cap f^{-1}(V\setminus\widehat\Delta)
$
such that $x_n\to x$, $f(x)=0$, and
\[
\ker d_{x_n}(f|_{S_\beta})\subset T_{x_n}\s_\epsilon,
\]
where $S_\beta$ is the stratum containing $x_n$. Passing to a subsequence, we may assume that all $x_n$ lie in the same stratum $S_\beta$ and that
$
\ker d_{x_n}(f|_{S_\beta})
$
converges in the Grassmannian to a subspace $T$. Then $T\subset T_x\s_\epsilon$. By Thom's condition,
\[
T_xS_\alpha=\ker d_x(f|_{S_\alpha})\subset T\subset T_x\s_\epsilon,
\]
where $S_\alpha\subset f^{-1}(0)$ is the stratum containing $x$, contradicting the choice of $\epsilon$.
\end{proof}

When $k=1$, given a holomorphic function $f\colon(X,\0)\to(\C,0)$, we can always choose a 
small enough representative such that $\widehat\Delta$ is either $\{0\}$ or  empty for any stratification. 
Moreover, we know that any Whitney stratification satisfies Thom's condition (see \cite{BMM}), 
so  we recover L\^e’s theorem:

\begin{corollary}[\cite{L0}]\label{L0}
Any holomorphic function $f\colon(X,\0)\to(\C,0)$ admits a Milnor-L\^{e} fibration.
\end{corollary}

Another well known case where $f$ admits a Milnor-L\^{e} fibration is when 
$X$ is smooth and $f$ defines an isolated complete intersection singularity (ICIS):

\begin{corollary}[Hamm \cite{Hamm}] Any holomorphic mapping $f\colon(\C^{n+k},\0)\to(\C^k,0)$ 
which defines an ICIS admits a Milnor-L\^{e} fibration where $\widehat\Delta=f(C)$ and $C$ is the subset of critical points of $f$.
\end{corollary}

However, when $k>1$, there are holomorphic mappings $f\colon(X,\0)\to(\C^k,0)$ 
which do not admit  a Milnor-L\^{e} fibration, even if $X$ is smooth and $f$ 
a complete intersection (with non-isolated singularity). The following example is due to L\^{e} and it can 
be found in \cite{Sabbah}. 

\begin{example}\label{ex:Le} {\rm Consider $\Phi\colon(\C^3,0)\to(\C^2,0)$ given by
$\Phi(x,y,z)=(y,x^2-y^2z)$. Then $f$ does not admit a Milnor-L\^{e} fibration and hence, 
it does not admit a stratification which satisfies Thom's condition.}
\end{example}

\subsection{Functions with isolated singularities relatively to $f$}\label{sub:Thom}

Let $f\colon (X,\0)\to(\C,0)$ be holomorphic. By \cite{BMM}, $f$ satisfies the Thom condition, hence admits a Milnor-L\^{e} fibration \cite{L0}. We study functions $g\colon (X,\0)\to(\C,0)$ for which the polar curve is well defined and
$
\Phi=(g,f)\colon (X,\0)\to(\C^2,0)
$
also admits a Milnor-L\^{e} fibration.

Fix:

\begin{enumerate}
\item a representative $\Phi=(g,f)\colon X\to V$, where $V$ is an open neighborhood of $0$ in $\C^2$;

\item a Whitney stratification of $f\colon X\to\C$ with strata $\mathcal S=(S_\alpha)$ on $X$, and two strata $\{\C\setminus\{0\},\{0\}\}$ on $\C$.
\end{enumerate}

Let $C(\Phi)$ denote the critical locus of $\Phi$ with respect to $\mathcal S$.

\begin{definition}\label{defn:nondepraved}{\rm
We say that $g\colon (X,\0)\to(\C,0)$ has an \emph{isolated critical point} with respect to $f$ if $C(g)\cap f^{-1}(0)=\{0\}$, where  $C(g)$ is the critical locus of $g$ with respect to $\mathcal S$.
}
\end{definition}

Thus ``with respect to $f$'' refers here to the fixed Whitney
stratification adapted to $f$; it does not mean that $g$ is assumed
to have an isolated critical point on each fiber of $f$.

The main result of this section is that if $g$ has an isolated critical point with respect to $f$, then $\Phi=(g,f)$ satisfies the Thom condition and therefore admits a Milnor-L\^{e} fibration.

\begin{theorem}\label{thm:Thom}
Suppose that $g\colon (X,\0)\to(\C,0)$ has an isolated critical point with respect to $f\colon (X,\0)\to(\C,0)$. Then:

\begin{enumerate}
\item $\Gamma:=\overline{C(\Phi)\setminus f^{-1}(0)}$ is empty or a curve in $X$;

\item $\Phi|_\Gamma\colon \Gamma\to V$ is finite. In particular,
$
\Delta:=\Phi(\Gamma)
$
is an analytic curve in $V$;

\item $\Phi\colon X\to V$ admits a stratification $\{\mathcal S',\mathcal T\}$ such that $\mathcal S'$ refines $\mathcal S$, has discriminant \\
$
\widehat\Delta=\Delta\cup(\C\times\{0\}),
$
and satisfies the Thom condition.
\end{enumerate}

In particular, $\Phi=(g,f)\colon (X,\0)\to(\C^2,0)$ admits a Milnor-L\^{e} fibration.
\end{theorem}

\begin{proof} We assume that $U\subset\C^N$ is an open neighbourhood of $\0$ such that $X$ is closed analytic in $U$ and $\tilde g\colon U\to\C$ is a holomorphic extension of $g\colon X\to\C$. 
Since $g$ has isolated critical point respect to $f$, we also assume $C(g)\cap f^{-1}(0)=\{0\}$.

We consider the relative conormal bundle of $f\colon X\to\C$, defined as
\[
T_{X,f}^*U:=\bigcup T_{S_\alpha,f}^*U,
\]
where
\[
T_{S_\alpha,f}^*U:=\{(x,\xi)\in T^*U\ |\ x\in S_\alpha \text{ and }\xi\left(\ker d_x(f|_{S_\alpha})\right)=0\}.
\]
By  \cite{BMM}, the Thom condition holds for $f\colon X\to \C$ and this implies that $T_{X,f}^*U$ is closed analytic in $T^*U$. Let $d\tilde g\colon U\to T^*U$ be the differential mapping given by $x\mapsto(x,d_x\tilde g)$. It follows that the set
\[
C:=(d\tilde g)^{-1}(T_{X,f}^*U)
\]
is also closed analytic in $X$. 

Observe that given a point $x\in S_\alpha$, $x\in C$ if and only if $x$ is a critical point of the restriction of $g$ to the fibre $S_\alpha\cap f^{-1}(f(x))$. In particular, when $f(x)=0$, $S_\alpha\subset f^{-1}(0)$ and hence $x\in C$ if and only if $x$ is a critical point of the restriction of $g$ to $S_\alpha$. We have $C\cap f^{-1}(0)=C(g)\cap f^{-1}(0)=\{\0\}$ and thus $\dim C\le 1$.  

Otherwise, if $f(x)\ne0$, then $f|_{S_\alpha}$ is a submersion and hence, $x\in C$ if and only if $x$ is a critical point of $\Phi|_{S_\alpha}$. We deduce that $\Gamma=\overline{C\setminus f^{-1}(0)}=\overline{C\setminus \{\0\}}$ and also $\dim \Gamma\le 1$. Moreover, $T_{S_\alpha,f}^*U$ has codimension $N-1$ in $T^*U$, so $(d\tilde g)^{-1}(T_{S_\alpha,f}^*U)$ has also codimension $\le N-1$ in $U$. Thus, $\Gamma$ has dimension 1 (if not empty). This shows item 1. But item 2 also follows, for $\Phi^{-1}(0)\cap \Gamma\subset C\cap f^{-1}(0)=\{\0\}$, hence $\Phi|_{\Gamma}\colon\Gamma\to V$ is finite after shrinking the representatives, if necessary.

The proof of  item 3 is just an adaptation of the proof given in \cite[Theorem 2.6]{GLN} (see also \cite[Theorem 5.6.3]{thesis-Roberto}) in the particular case that $g$ is the restriction of a submersion. In the target $V\subset \C^2$ we consider the stratification $\mathcal T$ with four strata: $V\setminus \widehat\Delta$, $\Delta\setminus\{\0\}$, $\C\times\{0\}\setminus\{\0\}$ and $\{\0\}$. In the source $X$, we consider the stratification $\mathcal S'$ with strata of the form $S_\alpha \cap \Phi^{-1}(T)\setminus \Gamma$ and $S_\alpha \cap \Phi^{-1}(T)\cap \Gamma$, for each $S_\alpha\in \mathcal S$ and $T\in \mathcal T$.

The strata $T\in\mathcal T$ are smooth, after shrinking the representatives, if necessary.
We prove that all the strata $A\in\mathcal S'$ are also smooth. 
Consider first a stratum of the form $A=S_\alpha \cap \Phi^{-1}(T)\setminus \Gamma$, with $S_\alpha\in \mathcal S$ and $T\in \mathcal T$. We  have three subcases: (a) $A=\{\0\}$, which is smooth trivially, (b) $A\subset f^{-1}(0)\setminus\{0\}$, with $T=\C\times\{0\}\setminus\{\0\}$, so $A=S_\alpha$ is also smooth or (c) $A\subset X\setminus f^{-1}(0)$, in this case $T=V\setminus \widehat\Delta$, $A=S_\alpha\setminus\Gamma$ is smooth again. 
Otherwise, we consider a stratum $A=S_\alpha \cap \Phi^{-1}(T)\cap \Gamma$, which is either the point $\{\0\}$ or has dimension 1 and its closure is analyitic. By the curve selection lemma, we can take a smaller representative such that $A$ is smooth.

By construction, $\Phi$ maps strata of $\mathcal S'$ onto strata of $\mathcal T$. We have to show that $\Phi$ maps these strata submersively.  If $A=S_\alpha\cap\Phi^{-1}(T)\setminus\Gamma$, the only non-trivial case is when $S_\alpha\subset f^{-1}(0)\setminus\{\0\}$ and $T=\C\times\{0\}\setminus\{\0\}$. We see that $A=S_\alpha$ and $\Phi\colon S_\alpha\to T$ is a submersion, because $g|_{S_\alpha}$ has no critical points. Otherwise, $A=S_\alpha\cap\Phi^{-1}(T)\cap\Gamma$ is mapped by $\Phi$ onto $T$. Now, the only non-trivial case is when $T=\Delta\setminus\{\0\}$. But $A$ and $T$ have dimension 1 and $\Phi$ is finite on $\Gamma$, so we can assume that  $\Phi\colon A\to T$ is a local diffeomorphism.

It is clear that $\mathcal T$ is a Whitney stratification of $V$, but we have to show that $\mathcal S'$ satisfies the Whitney conditions. The case of a pair strata contained in $X\setminus \Gamma$ follows from \cite[Lemma 2.5]{GLN}. The case of a pair of strata in $\Gamma$ is trivial, since one of them must be the stratum $\{\0\}$. Hence, it only remains to consider the case of a pair of strata $A,B\in\mathcal S'$ such that $A\subseteq \overline{B}$, $B\subset X\setminus\Gamma$ and $A\subset\Gamma$ has dimension 1. However, in this case, the set of points $x$ in $A$ such that $B$ is not Whitney regular over $A$ at $x$ is analytic and has dimension 0 (see e.g. \cite[Proposition 2.6]{GWPL}). After shrinking the representatives if necessary, we can assume that $B$ is Whitney regular over $A$.

Finally, it only remains to check the Thom property for $\{\mathcal S',\mathcal T\}$. Let $A,B\in\mathcal S'$ strata such that $A\subseteq \overline{B}$. If $A\subset\Gamma$, the restriction $\Phi:A\to T$ is a local diffeomorphism and the Thom condition is satisfied automatically. Otherwise, suppose that $A\subset X\setminus\Gamma$, we must have also $B\subset X\setminus\Gamma$ and 
\begin{equation}\label{eqn:AB}
A=S_\alpha \cap \Phi^{-1}(T)\setminus \Gamma, \quad B=S_\beta \cap \Phi^{-1}(T')\setminus \Gamma,
\end{equation}
for some $S_\alpha,S_\beta\in\mathcal S$ and $T,T'\in\mathcal T$. 

We take a sequence $\{x_n\}$ in $B$ which converges to a point $x\in A$ such that $\ker d_{x_n}(\Phi|_{B})$ converges to some space $P$. We must show $\ker d_{x}(\Phi|_{A})\subset P$. It follows from \eqref{eqn:AB} that
\[
\ker d_{x}(\Phi|_{A})=\ker d_{x}(\Phi|_{S_\alpha})=\ker d_{x}(f|_{S_\alpha})\cap\ker d_x\tilde g
\]
and analogously,
\[
\ker d_{x_n}(\Phi|_{B})=\ker d_{x_n}(\Phi|_{S_\beta})=\ker d_{x_n}(f|_{S_\beta})\cap \ker d_{x_n}\tilde g.
\]
After taking a subsequence if necessary, we can assume that $\ker d_{x_n}(f|_{S_\beta})$ converges to some space $Q$. Since $f\colon X\to\C$ satisfies the Thom condition for $\mathcal S$, $\ker d_{x}(f|_{S_\alpha})\subset Q$. Moreover, $x\notin\Gamma=C$, so $d_x\tilde g(\ker d_{x}(f|_{S_\alpha}))\ne0$ and hence, $d_x\tilde g(Q)\ne0$. But this implies that
\[
P=\lim_n\left(\ker d_{x_n}(f|_{S_\beta})\cap \ker d_{x_n}\tilde g\right)=Q\cap \ker d_x(\tilde g)\supset \ker d_{x}(f|_{S_\alpha})\cap\ker d_x\tilde g,
\]
as required.
\end{proof}

\begin{definition}\label{SIS}{\rm
Suppose that $g\colon (X,\0)\to(\C,0)$ has an isolated critical point at $\0$ with respect to $f$. The curve $\Gamma$ in Theorem \ref{thm:Thom} is called the \emph{polar curve relative to $g$}, and its image
$
\Delta:=\Phi(\Gamma)
$
is called the \emph{Cerf diagram}.
}
\end{definition}

\begin{example}{\rm
Let $f\colon (X,\0)\to(\C,0)$ be as above, and let $\ell$ be a general linear form as in \ref{general forms}. Then $\ell\colon (X,\0)\to(\C,0)$ has an isolated critical point with respect to $f$.
}
\end{example}

\begin{example}{\rm
Let $X$ be smooth at $\0$ and assume $(g,f)$ defines an ICIS. Then $g$ has an isolated singularity with respect to $f$.
}
\end{example}

\begin{example}{\rm
Consider the mapping of Example \ref{ex:Le},
$
(g,f)\colon(\C^3,\0)\to(\C^2,0),
$
with $g(x,y,z)=y$ and $f(x,y,z)=x^2-y^2z$. Take the stratification of $f\colon\C^3\to\C$ with strata
$
S_1=\{x=y=0\}, \quad
S_2=\{x^2-y^2z=0\}\setminus S_1, \quad
S_3=\C^3\setminus(S_1\cup S_2).
$
Since $g|_{S_1}=0$, the function $g$ does not have an isolated critical point with respect to $f$.
}
\end{example}

\subsection{Morsifications and the theorem}

We introduce the notion of Morsification of $g\colon(X,\0)\to(\C,0)$ relative to $f\colon(X,\0)\to(\C,0)$. This is a relative version of a Morsification with respect to a stratified space $X$, as defined in \cite[Section 4]{BR}.

\begin{definition}\label{def:Morsification} {\rm Let $f\colon(X,\0)\to(\C,0)$ be holomorphic and choose any stratification of $f$ as in Subsection \ref{sub:Thom}.
Given another holomorphic function $g\colon(X,\0)\to(\C,0)$, a \emph{Morsification of $g$ relative to $f$} is a 1-parameter deformation
$G\colon(X\times\C,\0)\to(\C,0)$  such that if we set  $G(x,t)=g_t(x)$, then 
$g_0=g$ and for each Milnor fibre $F_f=\B_\epsilon\cap X\cap f^{-1}(u)$, there
there exists $\eta$, with $0<\eta\ll \delta\ll\epsilon\ll 1$ such  if $0<|t|<\eta$, $g_t$ is a 
stratified Morse function on $F_f$.}
\end{definition}

The following lemma can be proved by using similar arguments to those given in \cite[Proposition 4.2]{BR} for the existence of a Morsification of $g$ with respect to $X$.

\begin{lemma}\label{Morsification}
Suppose that $g\colon (X,\0)\to (\C,0)$ has an isolated critical point with respect 
to $f\colon (X,\0)\to(\C,0)$. Then there exists a Morsification of $g$ relative to $f$.
\end{lemma}

%
%
%
%
%

Let $g\colon (X,\0)\to(\C,0)$ have an isolated critical point relative to $f\colon (X,\0)\to(\C,0)$, and denote by $F_f$ and $F_{g,f}$ the Milnor fibers of $f$ and $\Phi=(g,f)\colon (X,\0)\to(\C^2,0)$, respectively. Let $g_t$ be a Morsification of $g$ relative to $f$ as in Definition \ref{def:Morsification}. For each stratum $S_\a$, let $\delta_\a$ be the number of critical points of $g_t$ on $S_\a\cap F_f$.

We now determine the topology of the Milnor fibre $F_{f}$ using the Milnor fiber $F_{g,f}$ 
and the restriction to $F_f$ of the Morsification of $g$.

\begin{theorem}\label{main2} 
  Let $X$ be a closed analytic subset of a neighborhood of $0$ in $\C^N$, and let $f\colon (X,\0)\to(\C,0)$ be a holomorphic map germ. Let $g\colon (X,\0)\to(\C,0)$ have an isolated critical point relative to $f$, and denote by $F_f$ and $F_{g,f}$ the fibers of $f$ and $(g,f)\colon (X,\0)\to(\C^2,0)$, respectively. 
Then $F_f$ is homeomorphic to the space obtained from $F_{g,f}\times \D^2$ by attaching, along each stratum $S_\a\cap F_f$, $\delta_\a$ copies of the Morse data
 $$(D_\alpha,\partial D_\alpha)\times D_\alpha \times (\rm {Cone} \,(L_\a^+ \cup_{\partial} L_\a^-), L_\a^-) \,,$$ 
 where
 $\delta_\alpha$ is the number of Morse critical points of  a Morsification  $g_t$ of $g$ with respect to $f$, 
 
 $L_\a^+ \cup_{\partial} L_\a^-$ is homeomorphic to the link of the stratum $S_\a \cap F_f$, and $L^{\pm}_a$ are the half links;
 $D_\alpha$ is a real ball of dimension $n_\alpha=\dim (S_\alpha\cap F_f)$.  
 \end{theorem}

\begin{proof} 
We take a representative $\Phi=(g,f)\colon X\to V$.
By Theorem \ref{thm:Thom}, there exist $\epsilon,\delta$, with $0<\delta\ll\epsilon\ll1$, such that the restriction
\begin{equation}\label{ML2}
\Phi\colon \B_\epsilon \cap X\cap \Phi^{-1}(\rond{\D}_\delta\setminus \widehat
\Delta)\longmapsto \rond{\D}_\delta\setminus \widehat\Delta,
\end{equation}
is a locally trivial stratified fibration, with $\widehat\Delta=\Delta\cup(\C\times\{0\})$, and whose fiber is 
$$
F_{g,f}=\B_\epsilon \cap X\cap f^{-1}(z)\cap g^{-1}(u),
$$ 
for $(u,z)\in\rond{\B}_\delta\setminus \widehat\Delta$. 

We take $z\in\C$, with $0<|z|<\delta$ such that $z$ is a regular value of $f:
\B_\epsilon \cap X\to\C$. Then, there exists $\rho>0$ such that for all $u\in\C$, 
with $0<|u|<\rho$, $(u,z)\notin \widehat\Delta$. This implies that $(u,z)$ is a regular 
value of $\Phi$ and hence, $u$ is a regular value of $g\colon f^{-1}(z)\cap\B_\epsilon \cap X\to\C$.

Consider $G\colon (f^{-1}(z)\cap\B_\epsilon \cap X)\times\C\to\C$, given by $G(x,t)=g_t(x)$. 
We can choose now $\rho'$ with $0<\rho'<\rho$ such that for all $u$, with $0<|u|<\rho'$, $u$ is also a regular value of $G$.

Finally, let $\pi\colon G^{-1}(u)\cap \left((f^{-1}(z)\cap\B_\epsilon \cap X)\times\C\right)\to\C$ 
be the projection $\pi(x,t)=t$. Then $t$ is a regular value of $\pi$ if and only if $u$ is a regular 
value of $g_t\colon f^{-1}(z)\cap\B_\epsilon \cap X\to\C$. In particular, $0$ is a regular value of $\pi$ and hence, there exists $\eta>0$ such that 
$\pi\colon G^{-1}(u)\cap\left( (f^{-1}(z)\cap\B_\epsilon \cap X)\times \D_\eta\right)\to\D_\eta$ 
is a stratified proper submersion. By Thom-Mather's First Isotopy Lemma \cite[\S11]{Mather}, $\pi$ 
is a locally trivial fibration on $\D_\eta$. We deduce that for $0<|t|<\eta$, the fiber $F_{g,f}\equiv\pi^{-1}(0)$ is homeomorphic to 
$$
F_{g_t,f}=\B_\epsilon \cap X\cap f^{-1}(z)\cap g_t^{-1}(u)\equiv\pi^{-1}(t). 
$$
At this point, we can apply now the same arguments of stratified Morse theory as in the proof of 
Theorem \ref{main}, but using the function $g_t$ instead of the generic linear function $\ell$.
\end{proof}

\begin{corollary}
 We have:
\[
\chi(F_f)-\chi(F_{g,f})= \sum_{\alpha} \delta_\alpha m_\alpha \,,
\]
where   $m_\alpha=(-1)^{n_\alpha} \chi({\rm Cone}({\mathcal L}_\alpha),{\mathcal L}_\alpha)$,  and ${\mathcal L}_\alpha$ is the complex link of the stratum $S_\a \cap F_f$.
\end{corollary}

\begin{proof}
 We know from \cite {Sul} that the Euler characteristic is additive. Hence, 
Theorem \ref {main2}  implies:
$$\chi(F_f) = \chi(F_{g,f}) + \sum_{\alpha} \chi\Big      (D_\alpha,\partial D_\alpha)\times D_\alpha \times \big(\rm {Cone} \,(L_\a^+ \cup_{\partial} L_\a^-), L_\a^- \big)       \Big).
$$

At each stratum we have:
$$ \chi\Big   (   (D_\alpha,\partial D_\alpha)\times D_\alpha \times ({\rm {Cone}} \,(L_\a^+ \cup_{\partial} L_\a^-), L_\a^-) \Big) = (-1)^{n_\alpha} \; \chi
\big ({\rm {Cone}} \,(L_\a^+ \cup_{\partial} L_\a^-), L_\a^-)  \big)  \,
$$
and the result follows because for complex varieties the normal Morse data is homotopically equivalent to $(\rm Cone({\mathcal L}_\alpha),{\mathcal L}_\alpha)$, see for instance \cite[Corollary 1, p. 166]{GM}.
\end{proof}

\begin{example} {\rm Consider in $\C^n$ the variety $X=X_1\cup X_2$, where $X_1$ is the line $x_1=\dots=x_{n-1}=0$ and $X_2$ is the hyperplane $x_n=0$. We take two linear functions $f,g\colon X\to\C$ given by $f(x)=x_1+x_n$ and $g(x)=a_1x_1+\dots+a_nx_n$, where $a_1,\dots,a_n\in\C$ are generic.

The fibre of $f$ has two connected components $F_f=(F_f\cap X_1)\sqcup (F_f\cap X_2)\cong\{*\}\sqcup D^{2(n-2)}$, hence $\chi(F_f)=2$. But the fibre of $(g,f)$ has only one connected component  $F_{g,f}=F_{g,f}\cap X_2\cong D^{2(n-3)}$, which gives $\chi(F_{g,f})=1$ and thus, $\chi(F_f)-\chi(F_{g,f})=1$.

On the other hand $g$ has only one critical point on $F_f$ corresponding to the one single point stratum. It follows that $\delta_\alpha=1$ and $m_\alpha=1$, in concordance with our formula.
}
\end{example}

The following two examples are taken from \cite[Examples 4.1 and 4.2]{NP}. 

\begin{example}{\rm
Let $X$ be the surface in $\C^3$ parameterised as the image of the mapping $\C^2\to\C^3$ given by
\[
(t,s)\mapsto (t^3+st,t^4+(5/4)st^2,s).
\]
The implicit reduced equation of $X$ is:
\[
-4 s^4 y+5 s^3 x^2-32 s^2 y^2-16 s x^2 y+64 x^4-64 y^3=0,
\]
where $(x,y,s)$ are the coordinates in $\C^3$. We take the functions $f,g\colon X\to\C$ given by $f(x,y,s)=s$ and $g(x,y,s)=x$, so $\Phi=(g,f)\colon X\to\C^2$ 
 is the projection onto the $(x,s)$-plane.

The special fiber $f=0$ in $X$ is the plane curve $E_6$. The Milnor fiber $F_f=f^{-1}(s)\cap \mathbb B_{\varepsilon}$, with $0<|s|<\delta$, is the disentanglement, that is, the image of a stable perturbation of $E_6$. It has three nodes and hence the homotopy type of a bouquet of three $1$-spheres, so $\chi(F_f)=1-3=-2$ (see Figure \ref{fig:E6}). The mapping $\Phi=(g,f)\colon X\to\C^2$ is finite and has local degree three at the origin. Therefore, $F_{g,f}$ is equal to three points and $\chi(F_{g,f})=3$. We get $\chi(F_f)-\chi(F_{g,f})=-5$.

\begin{figure}[ht]
\begin{center}
\includegraphics[width=10cm]{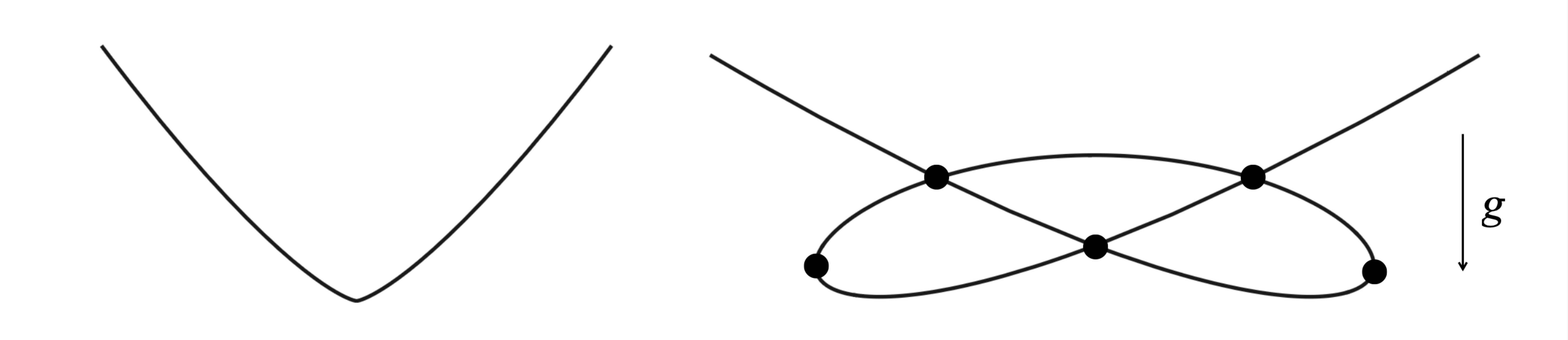}
\caption{}\label{fig:E6}
\end{center}
\end{figure}

To compute the sum $\sum_{\alpha} \delta_\alpha m_\alpha$, we see that $F_f$ has three single point strata which contribute with three critical points and other two critical points on the smooth part of $F_f$ (see again Figure \ref{fig:E6}). Each single point stratum is a node singularity and has complex link $\mathcal L_\alpha\simeq S^0$, which gives $m_\alpha=(-1)^0 \chi(Cone(\mathcal L_\alpha),\mathcal L_\alpha)=-1$. For each critical point on the smooth part we have $m_\alpha=(-1)^1 \chi(Cone(\mathcal L_\alpha),\mathcal L_\alpha)=-1$. This gives $\sum_{\alpha} \delta_\alpha m_\alpha=-5$.

}
\end{example}

\begin{example}{\rm 
Consider $X$ as the hypersurface in $\C^4$ parametrised as the image of the mapping $\C^3\to\C^4$:
\[
(u,v,w)\mapsto (u,v,w^2,w^5+uvw^3+(u^3-5uv-v)w).
\]
A simple computation shows that $X$ has the following reduced equation: 
\[
y^2-x(x^2+uvx+u^3-5uv-v)^2=0,
\]
 where $(u,v,x,y)$ are the coordinates in $\C^4$. Now we take $f\colon X\to \C$ 
 and $g\colon X\to\C$ as the functions $f(u,v,x,y)=v$ and $g(u,v,x,y)=u$, so $\Phi=(g,f)\colon X\to\C^2$ 
 is the projection onto the $(u,v)$-plane.
 
For $v\ne0$, the fibre $F_f=f^{-1}(v)\cap \mathbb B_{\varepsilon}$ is the disentanglement 
of the singularity $F_4$ in Mond's list \cite{Mond} 
(see Figure \ref{fig:F4}). Hence $F_f$ has the homotopy type of a bouquet of four $2$-spheres, so that $\chi(F_f)=1+4=5$. 

\begin{figure}[ht]
\begin{center}
\includegraphics[width=6cm]{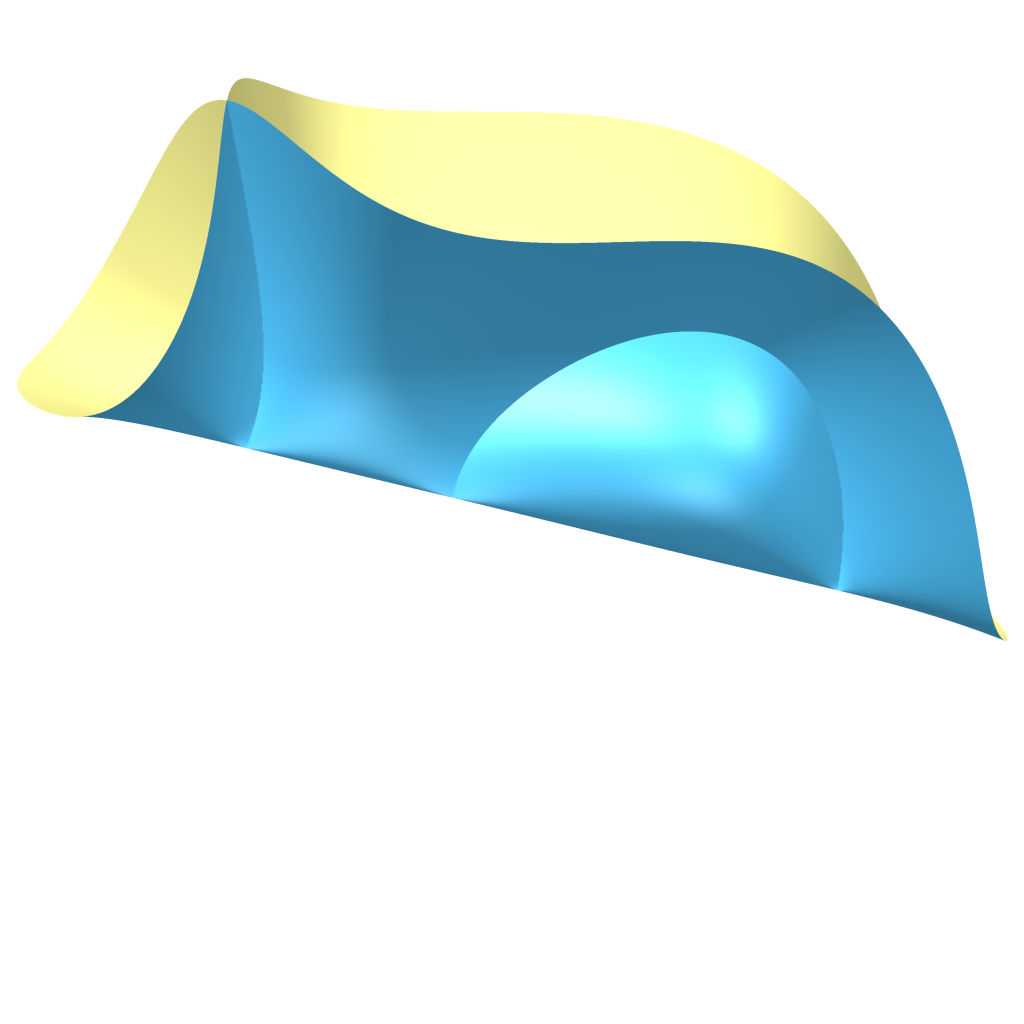}

\vskip-3cm\ 
\caption{}\label{fig:F4}
\end{center}
\end{figure}

On the other hand, for generic $(u,v)$, $F_{g,f}=g^{-1}(u)\cap f^{-1}(v)\cap \mathbb B_{\varepsilon}$ 
is the Morsification ot the plane curve singularity $A_4$, so it has the homotopy type of a bouquet 
of two $1$-spheres and hence, $\chi(F_{g,f})=1-2=-1$. We get $\chi(F_f)-\chi(F_{g,f})=6$.

Finally, by \cite[Example 4.2]{NP}, there are exactly six critical points of $g$ on each stratum in the interior of the fiber $F_f$. Three correspond to the $0$-dimensional stratum and arise from the Whitney umbrellas appearing in a stable perturbation of the $F_4$ singularity. The complex link $\mathcal L_\a$ of the Whitney umbrella has the homotopy type of $S^1$, hence
\[
m_\alpha=(-1)^{n_\alpha}\chi(\hbox{Cone} ({\mathcal L}_\alpha),{\mathcal L}_\alpha)=(-1)^0(1-0)=1.
\]
The other three critical points correspond to critical points of $g$ on the 1-dimensional stratum, 
which is the double point curve. At these points the complex link is $\mathcal L_\alpha\simeq S^0$, hence
\[
m_\alpha=(-1)^{n_\alpha}\chi(\hbox{Cone} ({\mathcal L}_\alpha),{\mathcal L}_\alpha)=(-1)^1(1-2)=1.
\]
Moreover, there are no critical points of $g$ on the 2-dimensional stratum of regular points. 
This gives  $\sum_{\alpha} \delta_\alpha m_\alpha=6$.
}
\end{example}

\section{An algebraic formula}\label{s: algebraic formula}
We now give an algebraic formula for the difference $\chi(F_f)-\chi(F_{g,f})$ in the case where  $\Sing(X)\subseteq f^{-1}(0)$, while $f$ may have non-isolated singularities. In particular, the formula holds when $X$ has isolated singularity. The formula is simple and has several applications, discussed in the next section.

Assume that $(X,\0)$ has pure dimension $d$, with $X\subset \C^N$. For a map germ $f\colon (X,\0)\to(\C^r,0)$, denote by $J_X(f)$ its Jacobian ideal with respect to $X$, defined as follows. Let $\phi_1,\dots,\phi_k$ generate the ideal $I_{X,0}$ in $\mathscr O_N$, and let $\bar f\colon (\C^N,\0)\to(\C^r,0)$ be a holomorphic extension of $f$. Then $J_X(f)$ is the ideal in $\mathscr O_{X,\0}$ given by
\[
J_X(f)=\frac{I_{X,\0}+J_{N-d+r}(\phi,\bar f)}{I_{X,\0}},
\]
where $J_{N-d+r}(\phi,\bar f)$ is the ideal in $\mathscr O_N$ generated by the $(N-d+r)$-
minors of the Jacobian matrix of $(\phi,\bar f)=(\phi_1,\dots,\phi_k,	\bar f_1,\dots,\bar f_r)$. 
When $(X,\0)$ is smooth, the zero locus $V(J_X(f))$ is precisely the set of 
critical points of $f\colon (X,\0)\to(\C^r,0)$.

Recall that given ideals $I,J$ in a ring $R$, the saturation of $I$ with respect to $J$ is the ideal
\[
I\colon J^\infty=\left\{a\in R\ |\ aJ^k\subseteq I, \text{ for some $k\ge 1$}\right\}.
\]
When $R=\mathscr O_{X,0}$ we have a nice geometrical interpretation: 
$$V(I\colon J^\infty)=\overline{V(I)\setminus V(J)}\,,$$
where the bar means the Zariski closure (see e.g. \cite{Singular}).

\begin{theorem}\label{thm:algebraic} Suppose that $(X,\0)$ has pure dimension $n+1$ and $\Sing(X)\subseteq f^{-1}(0)$ for some holomorphic function $f\colon (X,\0)\to(\C,0)$.
If $g\colon (X,\0)\to (\C,0)$ has an isolated critical point with respect to $f\colon (X,\0)\to(\C,0)$. Then,
\[
\chi(F_f)-\chi(F_{g,f})=(-1)^n \dim_\C\frac{\mathscr O_{X,\0}}{(f)+J_X(f,g)\colon f^\infty} \;.
\]
\end{theorem}

\begin{proof} Since $\Sing(X)\subseteq f^{-1}(0)$, we can assume that $X\setminus f^{-1}(0)$ is a stratum. 
The fibre $F_f$ is in this case a smooth manifold with boundary, 
hence it follows from Theorem \ref{main2} (see also \cite{CMSS}) that
\[
\chi(F_f)-\chi(F_{g,f})=(-1)^n \delta_\alpha,
\]
where $\delta_\alpha$ is the number of critical points of a Morsification of $g$ on 
the interior of $F_f$. By the conservation of the Milnor number, $\delta_\alpha$ 
is equal to the sum of the Milnor numbers of the critical points of $g$ on the interior 
of $F_f$. That is, for $0<\delta\ll\epsilon\ll 1$ and for $u\in\C$, $0<|u|<\delta$, 
the fibre $F_f=f^{-1}(u)\cap \B_\epsilon$ is smooth at $x\in \rond{F_f}=f^{-1}(u)\cap \rond{\B}_\epsilon$, thus:
\[
\delta_\alpha=\sum_{x\in \rond{F_f}} \mu(g|_{F_f};x)=\sum_{x\in \rond{F_f}}\dim_\C\frac{\mathscr O_{F_f,x}}{J_{F_f}(g)}.
\]
Moreover, $X$ is also smooth at any point $x\in \rond{F_f}$, which implies that
\[
\frac{\mathscr O_{F_f,x}}{J_{F_f}(g)}\cong \frac{\mathscr O_{X,x}}{(f-u)+J_X(f,g)}
\]

Suppose now that at $x=\0$ the ideal $J_X(f,g)$ has an irredundant primary decomposition in $\mathscr O_{X,\0}$ given by
\[
J_X(f,g)=\mathfrak q_1\cap\dots\cap\mathfrak q_t.
\]
On the one hand, the saturation with respect to $(f)$ can be written as
\[
J_X(f,g)\colon f^{\infty}=\bigcap_{i} (\mathfrak q_i\colon f^{\infty}).
\]
A simple computation shows that $\mathfrak q_i\colon f^{\infty}=\mathscr O_{X,\0}$ when $f\in \sqrt{\mathfrak q_i}$ and 
$\sqrt{\mathfrak q_i\colon f^{\infty}}=\sqrt{\mathfrak q_i}$ otherwise.
In particular, $J_X(f,g)\colon f^{\infty}$ does not have any $\mathfrak m$-primary 
component, that is, all its primary components have dimension $\ge 1$. 

On the other hand, we have
\[
V(J_X(f,g)\colon f^{\infty})=\overline{V(J_X(f,g))\setminus V(f))}=\Gamma,
\]
the relative polar curve. We consider in $\Gamma$ the analytic structure given 
by the ideal $J_X(f,g)\colon f^{\infty}$, so that $\Gamma$ has local ring 
$\mathscr O_\Gamma=\mathscr O_{X,\0}/J_X(f,g)\colon f^{\infty}$.

By item (1) in Theorem \ref{thm:Thom}, $\dim \Gamma=1$, so all primary components of $J_X(f,g):f^{\infty}$ have dimension $\le 1$. Hence $\mathscr O_\Gamma$ is unmixed of dimension $1$. 

Moreover, since $\Gamma\cap V(f)=\{\0\}$, the function $f$ defines a regular element in $\mathscr O_\Gamma$, and thus $\mathscr O_\Gamma$ is Cohen--Macaulay (see \cite[Exercise C.3.1]{Mond-NB}).

By \cite[Proposition 7.1]{Fulton}, we compute the local intersection number at $\0$ 
of the hypersurface $V(f)$ and the relative polar curve $\Gamma$ as the length:
\[
i(V(f),\Gamma;\0)=\dim_\C\frac{\mathscr O_{\Gamma}}{(f)}.
\]
Finally, we use  the conservation of the multiplicity (see \cite[Corollary E.6]{Mond-NB}):
 for $0<\delta\ll\epsilon\ll 1$ and for $u\in\C$, $0<|u|<\delta$,
\begin{align*}
i(V(f),\Gamma;\0)&=\sum_{x\in f^{-1}(u)\cap \rond{\B}_\epsilon} i(V(f-u),\Gamma;x)\\
&=\sum_{x\in f^{-1}(u)\cap \rond{\B}_\epsilon}\dim_\C\frac{\mathscr O_{\Gamma,x}}{(f-u)}\\
&=\sum_{x\in f^{-1}(u)\cap \rond{\B}_\epsilon}\dim_\C\frac{\mathscr O_{X,x}}{(f-u)+J_X(f,g)\colon f^\infty}\\
&=\sum_{x\in f^{-1}(u)\cap \rond{\B}_\epsilon}\dim_\C\frac{\mathscr O_{X,x}}{(f-u)+J_X(f,g)}=\delta_\alpha,
\end{align*}
since $f$ is a unit in the ring $\mathscr O_{X,x}$ when $x\in f^{-1}(u)\cap \rond{\B}_\epsilon$.
\end{proof}

\begin{remark}{\rm 
With the same hypothesis as in Theorem \ref{thm:algebraic}, Dutertre and Grulha in 
\cite[Corollary 5.4]{DN-arxiv}, obtain a similar formula:
\[
\chi(X\cap f^{-1}(\delta)\cap \B_\epsilon)-\chi(X^g\cap f^{-1}(\delta)\cap \B_\epsilon)=(-1)^{n} i(V(f),\Gamma;\0),
\]
where $X^g$ the hypersurface defined by the zeros of $g$ in $X$ and $0<|\delta|\ll\epsilon\ll1$. Thus, $X\cap f^{-1}(\delta)\cap \B_\epsilon$ is the Milnor fibre $F_f$, but $\chi(X^g\cap f^{-1}(\delta)\cap \B_\epsilon)$ is the Milnor fibre $F_{f|_{X^g}}$ of the restricted function $f|_{X^g}\colon (X^g,\0)\to(\C,0)$. They do not prove that the mapping $\Phi:=(g,f)\colon (X,\0)\to(\C^2,0)$ admits a Milnor-L\^e fibration and that in such case, the fibres $F_{g,f}$ and $F_{f|_{X^g}}$ coincide, although this is a consequence of our Theorem 
\ref{thm:Thom}. We remark that the above formula of \cite[Corollary 5.4]{DN-arxiv} appears only in the arXiv version of the paper, but not in the published version \cite{DN}.
}
\end{remark}

\begin{example}{\rm
Let $(g,f)\colon(\C^3,\0)\to(\C^2,\0)$ be the mapping $g(x,y,z)=x+y+z$ and $f(x,y,z)=xyz$. 
In this example, $f$ has a non-isolated critical point, but $g$ has an isolated  critical point
with respect to $f$. The fibre $F_f$ has the homotopy type of $\s^1\times\s^1$, so $\chi(F_f)=0$. 
The fibre $F_{g,f}$ is diffeomorphic to the Milnor fibre of the plane curve $xy(x+y)=0$ at the origin, 
so it has Milnor number 4 and  we have $\chi(F_{g,f})=1-4=-3$. Hence $\chi(F_f)-\chi(F_{g,f})=3$.

On the other hand, $J_2(f,g)$ is generated by the $2\times 2$-minors of the Jacobian matrix
\[
\left(
\begin{array}{ccc}
yz  & xz  & xy  \\
1  & 1  & 1  
\end{array}
\right)
\]
which are $yz-xz,yz-xy,xz-xy$. We get that $J_2(f,g)\colon f^\infty$ is generated by $y-z,y-x$ and hence
\[
\dim_\C\frac{\mathscr O_{3}}{(f)+J_2(f,g)\colon f^\infty}=\dim_\C\frac{\mathscr O_{3}}{(xyz,y-z,y-x)}=3\,.
\]
}
\end{example}

Suppose now that $X\subset\C^N$ has pure dimension $n+1$ and it has isolated 
singularity at $0$. Consider a holomorphic function $f\colon (X,\0)\to (\C,0)$ and  
choose generic linear functions $\ell_1,\dots,\ell_n:\C^N\to\C$. We define  $X_{n+1}=X$, $f_{n+1}=f$ and for $i=1,\dots,n$:
\[
X_i=X\cap\ell_n^{-1}(0)\cap\dots\cap\ell_i^{-1}(0),\quad f_i=f|_{X_i}\colon(X_i,\0)\to(\C,0).
\]

\begin{corollary}\label{cor:algebraic} With  the hypotheses of Theorem \ref{thm:algebraic}, one has:
\[
\chi(F_f)=\sum_{i=2}^{n+1} (-1)^{i-1} \dim_\C\frac{\mathscr O_{X_i,0}}{(f_i)+
J_{X_i}(f_i,\ell_{i-1})\colon f_i^\infty}+\dim_\C\frac{\mathscr O_{X_1,0}}{(f_1)}.
\]
\end{corollary}

\begin{proof} 
The genericity of the linear forms implies that for each $i=2,\dots,n+1$, $(X_i,\0)$ 
has pure dimension $i$ and an isolated singularity,  and $\ell_{i-1}\colon (X_i,\0)\to (\C,0)$ 
has an isolated critical  with respect to $f_i\colon (X_i,\0)\to(\C,0)$. Moreover, 
we have $F_{f_i,\ell_{i-1}}=F_{f_{i-1}}$. By Theorem \ref{thm:algebraic},
\[
\chi(F_{f_i})-\chi(F_{f_{i-1}})=(-1)^{i-1} \dim_\C\frac{\mathscr O_{X_i,\0}}{(f_i)+J_{X_i}(f_i,\ell_{i-1})\colon f_i^\infty}.
\]
Finally, for $i=1$, $(X_1,\0)$ has dimension 1 and $\chi(F_{f_1})$ is precisely the local 
degree of $f_1\colon(X_1,\0)\to(\C,0)$, that is,
\[
\chi(F_{f_1})=\deg(f_1)=\dim_\C\frac{\mathscr O_{X_1,\0}}{(f_1)}.
\]

\end{proof}

\begin{example} {\rm The formula given in Corollary \ref{cor:algebraic} can be 
easily implemented in \textsc{Singular} \cite{Singular2}. For instance, consider 
$X$ as the union of the 2-planes $x,y$ and $z,t$ in $\C^4$. We see that $X$ 
has pure dimension 2 and isolated singularity at $0$, although it is not Cohen-Macaulay. 
We also consider the function $f\colon(X,\0)\to(\C,0)$ given by $f(x,y,z,t)=(x+y+z+t)^3$. 
The following \textsc{Singular} code can by used to compute $\chi(F_f)$ with the formula of Corollary \ref{cor:algebraic}:

\texttt{
LIB ``primdec.lib'';\\ 
ring r=0,(x,y,z,t),ds;\\
ideal i1=x,y;\\
ideal i2=z,t;\\
ideal X=intersect(i1,i2);\\
poly f=(x+y+z+t)\^{}3;\\
poly g=3x-2y+z-t;\\
ideal F=X,f,g;\\
module M=jacob(F);\\
ideal j=minor(M,4)+X;\\
list l=sat(j,f);\\
ideal jj=l[1],f;\\
ideal Xg=X,g;\\
ideal X1=radical(Xg);\\
ideal F1=X1,f;\\
-vdim(std(jj))+vdim(std(F1));\\
6
}

The result is that $\chi(F_f)=6$. In fact, in this example is not difficult to see that $F_f$ 
is a disjoint union of six complex discs, which coincides with the computation.

}
\end{example}


\section{Applications}

In this section we discuss various consequences of Theorem \ref{thm:algebraic}.

\subsection{The classical L\^e-Greuel formula}\label{classical}

We follow the basic construction, as explained in Looijenga's book \cite[Section 5.A]{Looijenga}.
 We start with a holomorphic map germ $f=(f_1,\dots,k_k)\colon(\C^{n+k},\0)\to(\C^k,\0)$ 
 which defines an ICIS $(X,\0)=V(f_1,\dots,k_k)$ of dimension $n$ in $\C^{n+k}$. 
 After a linear change of coordinates in $\C^k$, we can assume that the $u_k$-axis in $\C^k$ 
 only meets the discriminant $\widehat\Delta$ of $f$ at $\0$. This generic choice of coordinates 
 implies that $f'=(f_1,\dots,f_{k-1})$ also defines an ICIS $(X',\0)=V(f_1,\dots,k_{k-1})$ of dimension $n+1$.
The classical L\^e-Greuel formula says:

\begin{theorem}[\cite{Greuel,L3}] With the above notation,
\[
\mu(X,\0)+\mu(X',\0)=\dim_\C\frac{\mathscr O_{n+k}}{(f_1,\dots,f_{k-1})+J_{k}(f_1,\dots,f_k)}.
\]
\end{theorem}

This formula follows easily from  Theorem \ref{thm:algebraic} as we show below. 

\begin{proof} We can perform a generic linear change of coordinates in $\C^{k-1}$ and 
assume that the $u_{k-1}$-axis only meets the discriminant $\Delta'_0$ of $f'$ at $\0$. 
Hence, $(X'',\0)=V(f_1,\dots,f_{k-2})$ is an ICIS of dimension $n+2$. But this choice of 
coordinates also forces that $F_f=F_{f''_{k-1},f''_{k}}$ and $F_{f'}=F_{f''_{k-1}}$ , 
where $f''_{k-1},f''_k\colon (X'',\0)\to(\C,0)$ are the restrictions of $f_{k-1},f_k$ to $(X'',\0)$, respectively. 

By Theorem \ref{thm:algebraic},
\[
\chi(F_{f'})-\chi(F_f)=(-1)^{n+1}\dim_\C \frac{\mathscr O_{X'',\0}}{(f''_{k-1})+J_{X''}(f''_{k-1},f''_{k})\colon (f''_{k-1})^\infty}.
\]
Since $f$ defines an ICIS, we have $C(f)\cap f^{-1}(\0)=\{\0\}$, where $C(f)$ is 
the critical locus of $f$. Hence, $C(f)$ is Cohen-Macaulay of dimension $k-1$ 
(see e.g. \cite[Proposition 4.7]{Looijenga}). Thus, the ideal $J_{X''}(f''_{k-1},f''_{k})$ 
defines a $1$-dimensional Cohen-Macaulay variety in $(X'',0)$. In particular, 
it is unmixed, that is, we have a primary decomposition
\[
J_{X''}(f''_{k-1},f''_{k})=\mathfrak q_1\cap\dots\cap \mathfrak q_r,
\]
where all its primary components $\mathfrak q_i$ have dimension 1. 
Also, $V(f''_{k-1})\cap V(J_{X''}(f''_{k-1},f''_{k}))=\{\0\}$, therefore $f''_{k-1}\notin\sqrt{\mathfrak q_i}$. We get:
\[
J_{X''}(f''_{k-1},f''_{k})\colon (f''_{k-1})^\infty=J_{X''}(f''_{k-1},f''_{k}),
\]
and hence
\begin{align*}
\chi(F_{f'})-\chi(F_f)&=(-1)^{n+1}\dim_\C \frac{\mathscr O_{X'',\0}}{(f''_{k-1})+J_{X''}(f''_{k-1},f''_{k})}\\
&=(-1)^{n+1}\dim_\C \frac{\mathscr O_{n+k}}{(f_1,\dots,f_{k-1})+J_{k}(f_1,\dots,f_k)},
\end{align*}
where the second equality follows easily from the definition of $J_{X''}(f''_{k-1},f''_{k})$. 
In this case, the fibres $F_{f'}$ and $F_f$ have the homotopy type of bouquet of spheres 
of dimensions $n$ and $n-1$, respectively, and the numbers of spheres in each case are 
the Milnor numbers $\mu(X',\0)$ and $\mu(X,\0)$, respectively. Thus, 
$F_{f'}=1+(-1)^n\mu(X',\0)$, $F_{f}=1+(-1)^{n-1}\mu(X,\0)$ and
\[
\chi(F_{f'})-\chi(F_f)=(-1)^{n+1}(\mu(X,\0)+\mu(X',\0)).
\]
\end{proof}

\subsection{L\^e-Greuel formula for smoothable space curves}\label{sub:curves}

In this subsection we assume that $(X,\0)$ is a space curve, that is, $\dim (X,\0)=1$ 
and $X\subset \C^N$ for some $N$. The Milnor number $\mu(X,\0)$ was introduced 
by Bassein \cite{Bassein} when $(X,\0)$ is smoothable and by Buchweitz-Greuel \cite{Buchweitz-Greuel} in general:
\[
\mu(X,\0)=\dim_\C\frac{\omega_{X,\0}}{d\mathscr O_{X,\0}},
\]
where $\omega_{X,\0}$ is the dualizing module of Grothendick and 
$d\colon \mathscr O_{X,\0}\to \omega_{X,\0}$ is induced by the exterior derivation.
When $(X,\0)$ is smoothable, the smoothing has the homotopy type of a bouquet of 
1-spheres and the number of such spheres is $\mu(X,\0)$.

Suppose we have a holomorphic function $f\colon(X,\0)\to(\C,0)$ which is finite, 
that is $f^{-1}(0)=\{\0\}$. The Milnor number $\mu(f)$ of a function on a space 
curve was introduced by Goryunov \cite{Goryunov} for space curves in $\C^3$ 
and by Mond-Van Straten \cite{Mond-VanStraten} in the general case. This is:
\[
\mu(f)=\dim_\C\frac{\omega_{X,\0}}{df(\mathscr O_{X,\0})},
\]
where now $df\colon \mathscr O_{X,\0}\to \omega_{X,\0}$ is induced by the exterior 
product with $df$. When $(X,\0)$ is smoothable, $\mu(f)$ is equal to the sum of the 
Milnor numbers of the critical points of $f$ on the smoothing. The main relationship between $\mu(f)$ and $\mu(X,\0)$ is:
\begin{equation}\label{eq:NT}
\mu(f)=\mu(X,\0)+\deg(f)-1,
\end{equation}
where $\deg(f)$ is the local degree, that is, the number of points of the fibre $F_f$ 
of $f\colon(X,\0)\to(\C,0)$ (see \cite{Nuno-Perez} when $(X,\0)$ is smoothable and \cite{Nuno-Tomazella} for the general case).

We assume that $(X,\0)$ is smoothable and take a smoothing, that is, a flat 
morphism $\pi\colon (\mathcal X,\0)\to(\C,0)$ such that $X_t:=\pi^{-1}(t)$ is 
smooth when $t\ne0$ and $X_0=X$. We denote by $\bar{f}:(\mathcal X,\0)\to(\C,0)$ 
any holomorphic extension of $f:(X,\0)\to(\C,0)$.
We give another proof of \eqref{eq:NT} based on our Theorems \ref{main2} and \ref{thm:algebraic}.

\begin{corollary}\label{cor:curves} With the above notation, we have:
\begin{align*}
\mu(f)&=\dim_\C\frac{\mathscr O_{\mathcal X,\0}}{(\pi)+J_{\mathcal X}(\pi,\bar{f})\colon \pi^\infty} \,,\\
\mu(f)&=\mu(X,\0)+\deg(f)-1 \,,\\
\mu(X,\0)&=\dim_\C\frac{\mathscr O_{\mathcal X,\0}}{(\pi)+J_{\mathcal X}(\pi,\ell)\colon \pi^\infty}-m_0(X,\0)+1 \,,
\end{align*}
where $\ell\colon\C^N\to\C$ is a generic linear form and $m_0(X,\0)$ is the multiplicity of $(X,\0)$. 
\end{corollary}

\begin{proof} The flatness of $\pi\colon (\mathcal X,\0)\to(\C,0)$ implies that 
$\dim(\mathcal X,\0)=2$ and that $\pi$ is regular in $\mathscr O_{\mathcal X,\0}$. 
But $f$ is finite, so it is also regular in $\mathscr O_{X,\0}=\mathscr O_{\mathcal X,\0}/(\pi)$. 
Thus $\pi,\bar{f}$ form a regular sequence and hence,  $\mathscr O_{\mathcal X,\0}$ is 
Cohen-Macaulay. In particular, $(\mathcal X,\0)$ has pure dimension 2. Moreover, the facts 
that $X_t$ is smooth for $t\ne0$ and $(X,\0)$ has isolated singularity give that $(\mathcal X,\0)$ has also isolated singularity.

We use Theorems \ref{main2} and \ref{thm:algebraic} with the functions $\pi,\bar{f}\colon(\mathcal X,\0)\to(\C,0)$:
\[
\chi(F_\pi)-\chi(F_{\bar f,\pi})=-\delta_\alpha=-
\dim_\C\frac{\mathscr O_{\mathcal X,\0}}{(\pi)+J_{\mathcal X}(\pi,\bar{f})\colon \pi^\infty} 
\]
where $\delta_\alpha$ is the number of critical points of a Morsification 
of $\bar{f}$ on the interior of $F_\pi$. But since $(X,\0)$ is smoothable, 
the number $\delta_\alpha$ is equal to $\mu(f)$ and we get the first formula. 

For the second formula, just observe that $\mu(X,\0)=1-\chi(F_\pi)$ and $\deg(f)=\chi(F_{\bar f,\pi})$.

Finally, the last formula follows from the other two when $f=\ell$, 
taking into account that the local degree of $\ell:(X,\0)\to(\C,0)$ is the multiplicity $m_0(X,\0)$.
\end{proof}

\subsection{L\^e-Greuel formula for isolated determinantal singularities}

We recall some basic definitions and properties of isolated determinantal 
singularities (IDS) from \cite{NOT}. Let $0<s\le m\le n $ be integers. We 
denote by $M_{m,n}=M_{m,n}(\C)$ the space of complex matrices of size 
$m\times n$ and by $M_{m,s}^s$ the subset of matrices with rank $<s$. 
We say that $(X,\0)$ is an IDS of type $(m,n;s)$ if $(X,\0)=F^{-1}(M_{m,n}^s)$ 
for some holomorphic map germ $F\colon(\C^N,\0)\to M_{m,n}$ such that:
\begin{enumerate}
\item $\codim (X,\0)=(m-s+1)(n-s+1)$,
\item $X$ is smooth at $x$ and $\rank F(x)=s-1$, for all $x\ne\0$ in a neighbourhood of the origin.
\end{enumerate}
For technical reasons we restrict ourselves to the cases where either $s=1$ 
or $N<(m-s+2)(n-s+2)$. Then, we have that $(X,\0)$ is smoothable. In fact, 
a determinantal smoothing is constructed as follows: choose a generic matrix $A\in M_{m,n}$ 
and consider the map germ $F_A\colon (\C^N\times\C,\0)\to M_{m,n}$ given by $F_A(x,t)=F(x)+tA$. 
Then, $(\mathcal X,\0)=F_A^{-1}(M_{m,n}^s)$ is a determinantal variety in $\C^N\times\C$ 
and $\pi\colon (\mathcal X,\0)\to(\C,0)$ is the restriction of the projection $\pi(x,t)=t$. 
If $A$ is generic, then for all $t\ne0$, $X_t:=\pi^{-1}(t)$ is smooth and $F_A$ has rank $s-1$ at any point in $X_t$. 

It follows that the topology of the Milnor fibre $F_\pi=X_t\cap \B_\epsilon$ of 
the determinantal smoothing is independent of the choice of $A$. In general, 
$F_\pi$ does not have the homotopy type of a bouquet of spheres (see \cite{Damon-Pike}), 
so we cannot speak of a Milnor number. We consider instead  the vanishing Euler characteristic 
defined in \cite{NOT} as $\nu(X,\0)=(-1)^d(\chi(F_\pi)-1)$, where $d=\dim (X,\0)$.

Assume now that $f\colon(X,\0)\to(\C,0)$ is a holomorphic function with isolated critical point 
defined on the IDS $(X,\0)$. We associate a pair of invariants with this function. As above, 
we denote by $\bar f\colon(\mathcal X,\0)\to(\C,0)$ any holomorphic extension of $f\colon(X,\0)\to(\C,0)$.
On the one hand, the Milnor number $\mu(f)$ is the number of critical points of $f_a|_{X_t}$, 
where $f_a$ is a Morsification of $\bar f$ on the determinantal smoothing $X_t$. 
On the other hand, the vanishing Euler characteristic of the fiber is also defined as
\[
\nu(X\cap f^{-1}(0),\0)=(-1)^{d-1}(\chi(F_\pi \cap f_a^{-1}(c))-1),
\]
where $c$ is a regular value of $f_a$. It is shown in \cite{NOT} that both numbers 
are well defined and depend only  on the germs of $X$ and $f$ at $\0$. 


In this setting we have again a L\^e-Greuel type formula for IDS. The following corollary 
can be proved by using the same arguments as in Corollary \ref{cor:curves}:

\begin{corollary}\label{cor:IDS} Let $f\colon(X,\0)\to(\C,0)$ be a holomorphic function 
with isolated critical point on an IDS $(X,\0)$ of dimension $d$. Then,
\begin{align*}
\mu(f)&=\dim_\C\frac{\mathscr O_{\mathcal X,\0}}{(\pi)+J_{\mathcal X}(\pi,\bar{f})\colon \pi^\infty},\\
\mu(f)&=\nu(X,\0)+\nu(X\cap f^{-1}(0),\0),\\
\nu(X,\0)&=(-1)^d\left(\sum_{i=2}^{d+1}(-1)^{i+1}\dim_\C\frac{\mathscr 
O_{\mathcal X_i,\0}}{(\pi_i)+J_{\mathcal X_i}(\pi_i,\ell_{i-1})\colon \pi_i^\infty}+m_0(X,\0)-1\right),
\end{align*}
where $\ell_1,\dots,\ell_d\colon\C^N\to\C$ are generic linear forms, 
$\mathcal X_i=\mathcal X\cap\ell_d^{-1}(0)\cap\dots\cap\ell_i^{-1}(0)$ and $\pi_i=\pi|_{\mathcal X_i}$.
\end{corollary}

 The second equation in Corollary \ref{cor:IDS} was proved originally in \cite{NOT}.

\subsection{L\^e-Greuel formula for image Milnor number}

The image Milnor number was introduced by Mond \cite{Mond2} for a 
hypersurface $(X,\0)$ in $\C^{n+1}$ which is defined as the image of 
a holomorphic map germ $f\colon(\C^n,S)\to(\C^{n+1},\0)$ with isolated instability. 
Here $S\subset \C^n$ is a finite subset whose cardinality is the number of irreducible components of $(X,\0)$ 
and the restriction $(\C^n,S)\to(X,\0)$ is the normalisation, since $f$ must be finite 
and generically one-to-one. The condition that $f$ has isolated instability ensures 
that $f$ admits a stabilisation, provided that either $f$ has corank one or $(n,n+1)$ 
are nice dimensions in the sense of Mather \cite{MatherVI}. By definition, a stabilisation 
is a 1-parameter unfolding $F\colon(\C^n\times\C,S)\to(\C^{n+1}\times\C,\0)$ given by 
$F(x,s)=(f_s(x),s)$ with the property that $f_0=f$ and $f_s$ has only stable singularities when $s\ne0$. 

Mond proved in \cite{Mond2} that $X_s\cap \B_\epsilon $ has the homotopy type 
of a bouquet of $n$-spheres, where $X_s$ is the image of $f_s$ and $0<|s|<\delta\ll\epsilon$. 
The number of such spheres is called the image Milnor number.  This  is denoted by $\mu_I(f)$ 
and it  is independent of the choice of $\delta,\epsilon$ and the stabilisation.  Since $\mu_I(f)$ 
is invariant under coordinate changes in the source and target of $f$, and the normalisation 
is unique up to coordinate changes in the source, it is also natural to consider it 
as an invariant of $(X,\0)$: the image Milnor number $\mu_I(X,\0)$ of $(X,\0)$.

There is a long standing conjecture by Mond \cite{Mond2} which claims that 
\[
\codim_{\mathscr A_e}(f)\le \mu_I(f),
\]
with equality when $f$ is weighted homogeneous. Here, $\codim_{\mathscr A_e}(f)$ 
is the extended $\mathscr A$-codimension and can be interpreted geometrically as 
the number of parameters in a miniversal unfolding of $f$. Hence, the conjecture is 
an inequality of type $\tau\le\mu$,   where $\tau$ is the usual Tjurina number of $(X,\0)$,
 but for the image Milnor number instead of the classical Milnor number. 
 The main difference is that, instead of deforming the implicit equation of $X$ 
 to get a smooth object, we deform the parametrisation of $X$ to get a stable object. 
 We refer to the survey \cite{NP} for a recent account on the Mond conjecture.

The image $(\mathcal X,\0)$ of the stabilisation $F$, together with the projection 
$\pi\colon(\mathcal X,\0)\to(\C,0)$, $\pi(y,s)=s$,  defines a flat deformation of 
$(X,\0)$ whose Milnor fibre is $F_\pi=X_s\cap \B_\epsilon$. We also have a natural 
stratification of $(\mathcal X,\0)$ given by the stable types on $\mathcal X\setminus\{\0\}$. 
We denote such stratification by $(S_\alpha)$. 

Assume now that $f$ has corank one and let $\ell\colon\C^{n+1}\to\C$ be a generic 
linear function, which defines a hyperplane $H=\ell^{-1}(0)$. Then 
$f^{-1}(H)$ is smooth at $S$ and the restriction $f\colon(f^{-1}(H),S)\to(H,\0)$ has 
isolated instability. As a consequence, its image $(X\cap H,\0)$ has also a well defined 
image Milnor number $\mu_I(X\cap H,\0)$. The following  L\^e-Greuel type formula in this context was proved in \cite{NP}:

\begin{theorem} With the above notation we have:
\[
\mu_I(X,\0)+\mu_I(X\cap H,\0)=\sum_{\alpha}\delta_\alpha,
\]
where $\delta_\alpha$ is the number of critical points of $\ell$ on the stratum $S_\alpha\cap X_s$.
\end{theorem}

The corank one assumption on $f$ also implies that $m_\alpha=1$ for every stratum $S_\alpha\cap X_s$ (see \cite{NP}). Thus, this formula is a special case of Theorem \ref{main2}. We do not yet have an algebraic version, as in the previous subsections. Such a version could help prove the Mond conjecture, at least in the corank one case.

\subsection{L\^e-Greuel formula for smoothable Gorenstein surfaces}

Let $(X,\0)$ be a normal surface singularity which is also smoothable and Gorenstein. Take a smoothing $\pi\colon(\mathcal X,\0)\to(\C,0)$ with Milnor fibre $F_\pi$. It follows 
from the works of Wahl \cite{Wahl} and Greuel and Steenbrink \cite{GS} 
that $\beta_1(F_\pi)=0$ and that $\beta_2(F_\pi)$ is independent of the choice of the smoothing, 
where $\beta_i$ is the $i$th-Betti number. Thus, it makes sense to call Milnor number to $\mu(X,\0)=\beta_2(F_\pi)$.

Suppose now that $f\colon(X,\0)\to(\C,0)$ is a holomorphic function with 
isolated critical point. Then the special fibre $(X\cap f^{-1}(0),\0)$ is 
a space curve and has a well defined Milnor number $\mu(X\cap f^{-1}(0),\0)$, as in Subsection \ref{sub:curves}. 
We choose any holomorphic extension $\bar{f}\colon(\mathcal X,\0)\to(\C,0)$ of $f$. We define 
the Milnor number $\mu(f)$ as the sum of the Milnor numbers of the critical points of $\bar{f}$ restricted 
to the interior of $F_\pi$.  

The next corollary follows again easily from Theorems \ref{main2} and \ref{thm:algebraic} 
and contains the L\^e-Greuel type formulas for smoothable Gorenstein surfaces:

\begin{corollary}\label{cor:Gorenstein} With the above notation we have
\begin{align*}
\mu(f)&=\dim_\C\frac{\mathscr O_{\mathcal X,0}}{(\pi)+J_\mathcal X(\pi,\bar{f})\colon \pi^\infty},\\
\mu(f)&=\mu(X,\0)+\mu(X\cap f^{-1}(0),\0),\\
\mu(X,\0)&=\dim_\C\frac{\mathscr O_{\mathcal X,0}}{(\pi)+
J_\mathcal X(\pi,\ell_2)\colon \pi^\infty}-\dim_\C\frac{\mathscr O_{\mathcal X_2,0}}{(\pi_2)
+J_{\mathcal X_2}(\pi_2,\ell_1)\colon \pi_2^\infty}+m_0(X,\0)-1,
\end{align*}
where $\ell_1,\ell_2:\C^N\to\C$ are generic linear forms, 
$\mathcal X_2=\mathcal X\cap\ell_2^{-1}(0)$ and $\pi_2=\pi|_{\mathcal X_2}$.
\end{corollary}

An immediate consequence of the corollary is that $\mu(f)$ is well defined and  depends only on $f\colon(X,\0)\to(\C,0)$.

\begin{example}{\rm 
We consider a surface $(X,\0)\subset\C^5$ defined by the Pfaffians $P_1,\dots,P_{2r+1}$ of a skew-symmetric matrix
$
F=(f_{ij})
$
of size $(2r+1)\times(2r+1)$, with entries $f_{ij}\in\mathscr O_5$ and $f_{ji}=-f_{ij}$. Recall that $P_i$ is the square root of the determinant of the submatrix obtained by deleting the $i$th row and column (see \cite{Buchsbaum-Eisenbud}).

If $(X,\0)$ has dimension 2, then it is Gorenstein, by a theorem due to 
Buchsbaum and Eisenbud \cite{Buchsbaum-Eisenbud}. In fact, every Gorenstein surface germ in $\C^5$ arises this way. When $r=1$, $(X,\0)$ is a complete intersection; if $r>1$ and all $f_{ij}$ lie in the maximal ideal $\mathfrak m_5$, then $(X,\0)$ is not. As far as we know, this is the simplest example of a Gorenstein surface that is not a complete intersection.

If the $f_{ij}$ are chosen so that $(X,\0)$ has an isolated singularity, then $(X,\0)$ is normal and smoothable. A smoothing is obtained by taking a generic skew-symmetric constant matrix
$
A=(a_{ij}),
$
with $a_{ji}=-a_{ij}$, and considering the $3$-fold $(\mathcal X,\0)\subset\C^5\times\C$ defined by the Pfaffians of $F+tA$, together with the projection 
$
\pi\colon(\mathcal X,\0)\to(\C,0),\; \pi(x,t)=t.
$
Then $\pi^{-1}(0)=X$, while $\pi^{-1}(t)$ is smooth for $t\neq0$ and $A$ sufficiently general.

We compute the Milnor number $\mu(X,\0)$ when $F$ is a $5\times5$ skew-symmetric matrix whose entries are linear forms. For generic coefficients, $(X,\0)$ has an isolated singularity. Using Corollary \ref{cor:Gorenstein} and 
\textsc{Singular}, we obtain
$
\mu(X,\0)=4
$ (see the code below).

\texttt{
LIB ``matrix.lib'';\\
LIB ``sing.lib'';\\
ring r=0,(x,y,z,w,v,t),ds;\\
ideal i = ideal(randommat(1,10,maxideal(1),9));\\
matrix A=skewmat(5,i);          // Skew-symmetric matrix of generic linear forms \\ 
poly Pf1=A[1,2]*A[3,4]-A[1,3]*A[2,4]+A[1,4]*A[2,3];   // The Pfaffians  \\
poly Pf2=A[1,3]*A[4,5]-A[1,4]*A[3,5]+A[1,5]*A[3,4];\\
poly Pf3=A[1,2]*A[4,5]-A[1,4]*A[2,5]+A[1,5]*A[2,4];\\
poly Pf4=A[1,2]*A[3,5]-A[1,3]*A[2,5]+A[1,5]*A[2,3];\\
poly Pf5=A[2,3]*A[4,5]-A[2,4]*A[3,5]+A[2,5]*A[3,4];\\
ideal X=Pf1,Pf2,Pf3,Pf4,Pf5;  // The ideal generated by the Pfaffians \\
dim(std(X));\\
3\\
mres(X,0);  // We check it is Gorenstein of codimension 3 \\
\hspace*{10pt}1\ \ \ \ \    5\ \ \ \ \        5\ \ \ \ \        1\\      
r <--  r <--  r <--  r\\
\\
0\ \ \ \ \        1\ \ \ \ \        2\ \ \ \ \        3 \\     
poly pi=t;\\
ideal X0=X,pi;\\
module M0=jacob(X0);\\
ideal J0=minor(M0,3)+X;\\
dim(std(J0)); // We check X0 has isolated singularity and that X is a smoothing\\
0\\
poly ell2=3x-2y+z-w+v;  // A generic linear form\\
ideal F=X,pi,ell2;\\
module M=jacob(F);\\
ideal J=minor(M,4)+X;\\
list L=sat(J,pi);\\
ideal K=L[1],pi;\\
ideal X2=X,ell2;\\
poly ell1=x+y-2z-3w+v; // Another generic linear form\\
ideal F2=X2,pi,ell1;\\
module M2=jacob(F2);\\
ideal J2=minor(M2,5)+X2;\\
list L2=sat(J2,pi);\\
ideal K2=L2[1],pi;\\
vdim(std(K))-vdim(std(K2))+mult(std(X0))-1; // The Milnor number\\
4\\
}
}
\end{example}



\subsection{A remark on indices of vector fields}
Let $(X,\0)\subset\C^N$ be a complex analytic germ of dimension $n+k$, and let
$f:(X,\0)\to(\C^k,0)$ be a holomorphic map-germ endowed with complex Whitney stratifications adapted to $V(f)$. Assume that $\{\0\}$ is a stratum and that $f$ satisfies the Thom condition. Then $f$ admits a local Milnor-L\^e fibration.

Let $B_\varepsilon$ be a Milnor ball for $X$, chosen so that every
stratum of $X\cap B_\varepsilon$ has $0$ in its closure, and let
\[
F_f=f^{-1}(t)\cap B_\varepsilon
\]
be a Milnor fiber of $f$.
Let $g:(X,0)\to(\mathbb C,0)$ be a germ with an isolated stratified
critical point relative to $V(f)$, and assume that $(f,g)$ satisfies
the Thom condition. Let $F_{f,g}$ be the corresponding Milnor fiber.
Both $F_f$ and $F_{f,g}$ inherit the stratification induced from $X$.

As in \cite[Section~2]{CMSS}, near $\partial F_f$ one projects the
gradient $\nabla g(x)$ onto the stratum containing $x$, and glues these
fields by a partition of unity to obtain a stratified vector field
$\nabla g|_{F_f}$. Extending it over the interior of $F_f$ by the
radial extension of M.~H.~Schwartz, let
$\operatorname{Ind}(\nabla g|_{F_f})$ be the sum of its local Schwartz
(radial) indices. The L\^e--Greuel formula of \cite{CMSS} states:
\[
\chi(F_f)-\chi(F_{f,g})
=
\operatorname{Ind}(\nabla g|_{F_f}).
\tag{7.2}
\]

It is useful to compare the hypotheses in \cite{CMSS} with those used
here. The result of \cite{CMSS} is formulated for maps
\[
f:(X,0)\longrightarrow(\mathbb C^k,0),
\]
with $k$ arbitrary, assuming that $f$ is generically a submersion with
respect to a Whitney stratification and satisfies the Thom
$a_f$-condition. The function $g$ is assumed to have an isolated
stratified critical point both on $X$ and on $V(f)$. In the present
paper we restrict to $k=1$ and start with a fixed Whitney
stratification adapted to $f$. We assume that $g$ has an isolated
critical point with respect to this stratification, without imposing
a separate isolated-critical-point assumption on $g|_{V(f)}$.
Theorem~5.7 then provides the geometric properties needed for the pair
\[
\Phi=(g,f):(X,0)\longrightarrow(\mathbb C^2,0):
\]
the relative critical locus away from $f^{-1}(0)$ is empty or a curve,
its image is an analytic discriminant curve, and $\Phi$ admits a
suitable stratification satisfying the Thom condition.

The problem is to compute
$\operatorname{Ind}(\nabla g|_{F_f})$ when $F_f$ is singular.
As observed in \cite{CMSS}, results of Dutertre--Grulha give a
stratified expression for this index in terms of the critical points
of a morsification of $g$ on the different strata and the topology of
the corresponding complex links. Thus, when $k=1$ and in the common
range of hypotheses, this gives
\[
\operatorname{Ind}(\nabla g|_{F_f})
=
\sum_\alpha \delta_\alpha m_\alpha,
\]
where $\delta_\alpha$ is the number of critical points on each stratum
of a Morsification of $g$, and
\[
m_\alpha=(-1)^{n_\alpha}
\chi(\operatorname{Cone}(L_\alpha),L_\alpha),
\]
with $L_\alpha$ the complex link of the stratum.

Consequently, in the common range of hypotheses, the
Euler-characteristic consequence of Theorem~5.13 agrees with the
previously known numerical L\^e--Greuel formula. The additional
content of Theorem~5.13 is topological: it describes $F_f$ as obtained
from $F_{g,f}\times D^2$ by attaching, along the different strata, the
corresponding stratified Morse data with multiplicities
$\delta_\alpha$. Thus the formula of \cite{CMSS} and Theorem~5.13
give complementary descriptions of the same Euler-characteristic
difference: the former identifies it with a generalized index, while
the latter exhibits the underlying topological modification of the
Milnor fiber.

Is this a special case of a more general theorem?


\end{document}